\documentclass[11 pt, psamsfonts]{amsart}
\usepackage{amssymb}

\def\ignore#1{\relax}

\def\g{\mathfrak g}
\def\h{\mathfrak h}

\def\sl{\mathfrak sl}
\def\so{\mathfrak so}

\def\Z{{\mathbb Z}}

\def\Q{{\mathbb Q}}
\def\C{{\mathbb C}}

\def\la{\lambda}

\def\r{{\bf r}}
\def\Ra{\mathcal R}
\def\N{\mathbb N}

\def\Ca{\mathcal C}

\def\Da{\mathcal D}

\def\b{{\bf b}}

\def\ignore#1{\relax}
\def\om{\omega}
\def\e{\epsilon}
\def\1{{\bf 1}}

\def\End{{\rm End}}

\calclayout

\renewenvironment{proof}[1][\proofname]{\par
  \normalfont
  \topsep6\p@\@plus6\p@ \trivlist
  \item[\hskip\labelsep\bfseries
    #1\@addpunct{.}]\ignorespaces
}{%
  \qed\endtrivlist
}

\def\th@plain{%
  \let\thmhead\thmhead@plain \let\swappedhead\swappedhead@plain
  \thm@preskip.5\baselineskip\@plus.2\baselineskip
                                    \@minus.2\baselineskip
  \thm@postskip\thm@preskip
  \itshape
\renewcommand{\labelenumi}{{(\alph{enumi})\quad}}
                        \renewcommand{\labelenumii}{{(\roman{enumii})\ }}
}
\def\th@definition{%
  \let\thmhead\thmhead@plain \let\swappedhead\swappedhead@plain
  \thm@preskip.5\baselineskip\@plus.2\baselineskip
                                    \@minus.2\baselineskip
  \thm@postskip\thm@preskip
  \upshape
}
\def\th@remark{%
  \thm@headfont{\itshape}
  \let\thmhead\thmhead@plain \let\swappedhead\swappedhead@plain
  \thm@preskip.5\baselineskip\@plus.2\baselineskip
                                    \@minus.2\baselineskip
  \thm@postskip\thm@preskip
  \upshape
}

{\theoremstyle{plain}
\newtheorem{theorem}{Theorem}[section]
}

{\theoremstyle{plain}
\newtheorem{proposition}[theorem]{Proposition}
}

{\theoremstyle{plain}
\newtheorem{corollary}[theorem]{Corollary}
}

{\theoremstyle{plain}
\newtheorem{lemma}[theorem]{Lemma}
}

{\theoremstyle{plain}

}

{\theoremstyle{definition}
\newtheorem{definition}[theorem]{Definition}
}

{\theoremstyle{definition}

}

{\theoremstyle{remark}
\newtheorem{remark}[theorem]{Remark}
}

{\theoremstyle{remark}

}

\numberwithin{equation}{section}

\renewcommand{\labelenumi}{{ \theenumi.}}
\renewcommand{\labelenumii}{{(\alph{enumii})}}

\def\la{\lambda}
\def\al{\alpha}

\def\choose #1 #2{\begin{pmatrix}#1\\#2\end{pmatrix}}

\def\Vnk{V_n^{(k)}}
\def\ek{e_{(k)}}
\def\WBk{W(B_k)}
\def\Zq{\Z[q,q^{-1}]}
\def\qB{Br}
\def\BnN{Br_n(N)}
\def\b{b}
\def\bb{\tilde b}
\def\gp{g^+}
\def\gm{g^-}

\def\eb{\bar e}
\def\ekb{\eb_{(k)}}
\def\Psik{\Psi_k}

\def\j{{\bf j}}
\def\Brf{{\tilde{Br}_4(r,q)}}
\def\qBb{\overline{\qB}}
\def\LNell{\Lambda(N,\ell)}
\def\LNellb{\bar\Lambda(N,\ell)}
\def\r{\rho}

\begin{document}

\title[A $q$-Brauer Algebra ]
{A $q$-Brauer Algebra }

\author{Hans Wenzl}
\thanks{}

\address{Department of Mathematics\\ University of California\\ San Diego,
California}

\email{hwenzl@ucsd.edu}

\begin{abstract}
We define a new $q$-deformation of Brauer's centralizer algebra
which contains Hecke algebras of type $A$ as unital subalgebras.
We determine its generic structure as well as the structure
of certain semisimple quotients. This is expected to have
applications for constructions of  subfactors of type II$_1$
factors and for module  categories of fusion categories
of type $A$ corresponding to certain symmetric spaces.
\end{abstract}
\maketitle

In his paper \cite{Brauer}, Richard Brauer introduced a series of algebras,
specializations of which describe the decomposition of tensor powers of
the defining vector representation of an orthogonal or
symplectic group. More recently, $q$-deformations of these algebras
have been defined in \cite{BW} and \cite{Mu} in connection with knot
theory and quantum groups. They found a number of applications,
such as in the study of subfactors and tensor categories
(see e.g. \cite{w1}, \cite{TW}, \cite{TuW}).

In this paper we introduce another $q$-deformation of Brauer's
centralizer algebras motivated by the following problem: Let $V$ be
the $N$-dimensional representation of $Gl(N)$. Restricting the
action of $Gl(N)$ on tensor powers $V^{\otimes n}$ to $O(N)$ leads
to embeddings of the centralizer algebras $\C S_n$, where $S_n$ is
the symmetric group,  into the Brauer algebra $D_n(N)$. Our idea now
is very simple: Find a $q$-deformation of $D_n(N)$ which extends the
$q$-deformation of $\C S_n$, the Hecke algebra $H_n(q)$ of type
$A_{n-1}$, subject to certain compatibility conditions with respect
to taking tensor products. This can also be stated in the language
of module categories (see the beginning of Section
\ref{motivation}). We shall see that these conditions completely
determine a $q$-deformation of the Brauer algebra $D_n(N)$. This
approach also carries over comparatively easily to the setting of
fusion tensor categories, i.e. for certain quotients of Hecke
algebras at roots of unity. This will be important for one of the
main motivations of this work, the constructions of examples of
subfactors of II$_1$ von Neumann factors. They were, at least in
part, inspired by work in conformal field theory in connection with
twisted affine loop groups and boundary conformal field theory
(see e.g. \cite{GG} and references
therein).

It is well-known that in this context the Hecke algebras correspond
to Jimbo-Drinfeld quantum groups $U_q\sl_N$
via an extension of
Schur duality. So our new algebras should correspond to
a $q$-deformation of the subalgebra $U\so_N\subset U\sl_N$.
Such algebras were defined as
coideal algebras  in work by Letzter (see \cite{L1}, \cite{L}),
and also in work by Gavrilik and Klimyk and by Noumi (see \cite{GK}, \cite{N}).
This could give another, potentially more conceptual
approach to derive our algebras, at least for the generic case
with $q$ not a root of unity. Related work in this direction
has already appeared earlier in \cite{Mo},
see the remarks at the end of this paper.
So our algebras can also be viewed as part of a categorical construction
of quantum analogs of certain symmetric pairs. Our approach also
works in the context of fusion categories, which, so far, would not be
so clear in the context of coideal algebras.

Here is a brief outline of the contents of this paper. In the first
section, we review results about Brauer's centralizer algebras
and Hecke algebras. This
will also serve as a model for our approach of defining and proving
results about our $q$-deformation of Brauer's centralizer algebra.
In the second section, we motivate our definitions via an approach
to find module categories of quantum groups from subalgebras
of the classical Lie algebra. We then
give the definition of our algebras via
generators and relations in the following section. We  show that they have
bases labeled by the basis graphs of Brauer's algebras. Moreover,
they also have the same decomposition into full matrix rings in
the generic case as Brauer's. In the fourth section, we define a trace
functional on our algebras with certain properties. It is an extension
of certain important trace functionals defined on Hecke algebras,
which are often referred to as Markov traces.
We will use our results on these Markov traces in Section 5 to
determine for which values
of the parameters our algebras will be semisimple. Moreover, we also
determine certain semisimple quotients in the non-semisimple case.
One can also see at these quotients that the algebras in this paper
are different from the $q$-deformations of Brauer algebras in \cite{BW}
and \cite{Mu}. We then discuss  several applications  of our algebras
such as the construction of module categories, subfactors and
 representations of fusion rings.

$Acknowledgements$: It is a pleasure to thank David Jordan,
Viktor Ostrik and Antony Wassermann for useful discussions and references.

\section{Brauer and Hecke algebras}

\subsection{Basic definitions}\label{basicdefinitions} In this paper
Brauer's centralizer algebra $D_n$ is defined over the ring $\Z[x]$
via  a basis given by  graphs with $2n$ vertices, arranged
at two levels, and $n$ edges, where each vertex belongs to exactly
one edge. We will call an edge $vertical$ if its vertices are
on different levels, and $horizontal$ if they are on the same level.
Concatenation of two basis graphs $a$ and $b$ is given similarly as
with braids. One puts $a$ on top of $b$ such that the $n$ lower vertices
coincide with the $n$ upper vertices of $b$. One then removes all
cycles, i.e. parts of the resulting graph which are not connected to
an upper or lower vertex. The element $ab$ is then defined to be
this resulting graph without cycles, multiplied by $x$ taken to the
power of the number of removed cycles; here $x$ is a variable.
To give an example, let $\ek$ be the element of $D_n$ given by
a graph with $2k$ horizontal edges on the left, and the remaining
$n-2k$ edges vertical. E.g. see below the graph for $e_{(2)}\in D_7$:
\vskip 0in
\begin{picture}(300,100)(0,0)
\put(140,50){\circle*{4}}
\put(140,50){\line(1,0){20}}
\put(160,50){\circle*{4}}
\put(140,70){\circle*{4}}
\put(140,70){\line(1,0){20}}
\put(160,70){\circle*{4}}
\put(180,50){\circle*{4}}
\put(180,50){\line(1,0){20}}
\put(200,50){\circle*{4}}
\put(180,70){\circle*{4}}
\put(180,70){\line(1,0){20}}
\put(200,70){\circle*{4}}
\put(220,50){\circle*{4}}
\put(220,50){\line(0,0){20}}
\put(220,70){\circle*{4}}
\put(240,50){\circle*{4}}
\put(240,50){\line(0,0){20}}
\put(240,70){\circle*{4}}
\put(260,50){\circle*{4}}
\put(260,50){\line(0,0){20}}
\put(260,70){\circle*{4}}
\put(180,10){Figure 1}
\end{picture}
\vskip 0in
Then it is easy to check that $\ek e_{(m)}=e_{(m)}\ek=x^m\ek$ for
any $m\leq k$; here the horizontal edges of $\ek$  should be drawn slightly concave
to obtain cycles.
In the following, Brauer's centralizer algebra $D_n$ is the free
$\Z[x]$-module spanned by the above mentioned basis graphs.
It is clear from the definition that the multiplication of $D_n$
is well-defined over $\Z[x]$ and associative. It is also clear
that its rank is $n!!=1\cdot 3\cdot\ ...\ (2n-1)$.

Observe that $D_n$ contains a subalgebra which is isomorphic to $\Z[x] S_n$,
where $S_n$ is the symmetric group of all permutations of $n$ symbols.
It is spanned by the basis graphs which only have vertical edges.
Then we get a decomposition
of $D_n(x)$ in terms of $S_n-S_n$ bimodules as
\begin{equation}\label{Snbimod}
D_n(x)\ \cong\ \bigoplus_{k=0}^{[n/2]} \Z[x] S_n\ek S_n;
\end{equation}
informally, $S_n\ek S_n$ can be viewed as the set of all graphs
with exactly $2k$ horizontal edges. Moreover, as the product of
two graphs has at least as many horizontal graphs as either of
them, it is easy to see that
$I(m)=\bigoplus_{k\geq m} \Z[x] S_n\ek S_n$
is a two-sided ideal in $D_n$ for each $m$ with $2m\leq n$.

It is clear from the pictures that multiplication of a graph
of $D_n$ from the left (i.e. from above pictorially) does not
change the position of the lower horizontal edges. This defines
a decomposition of $\Z[x] S_n\ek S_n$
into $S_n$-modules. Combinatorially, the position of the
lower horizontal edges of a graph in $ S_n\ek S_n$ is
determined as follows: We choose a subset of $2k$
elements from $[1,n]$ (only integers) and partition it into $k$ subsets
of 2 elements each. Let $P(n,k)$ be the set of all those partitions.
Then
\begin{equation}\label{Snmod}
 \Z[x] S_n\ek S_n\ \cong\ \bigoplus_{\j\in P(n,k)} \Z[x] S_n \ek w_\j,
\end{equation}
where $w_\j\in S_n$ such that $\ek w_\j$ is the graph whose lower
horizontal edges are given by the partition $\j\in P(n,k)$
and such that no vertical edges intersect. This completely determines
$\ek w_\j$. The permutation $w_\j$ is not uniquely determined. We shall
later make the choice of $w_\j$ more precise.

We shall also consider the Brauer algebra $D_n(N)$, $N\in\Z$ which
is defined over $\Z$ by the same graphs as before. The only difference
is that now the variable $x$ is replaced by the integer $N$.

\subsection{The module $\Vnk$ for Brauer algebras} It is also easy to
see that multiplication of a graph $d\in \Z[x] S_n \ek w_\j$ by an
element in $D_n$ from the left/above leaves the lower horizontal edges
unchanged, but may add additional lower horizontal edges.
Hence the factor module $\Z[x] S_n \ek w_\j +I(k+1)/I(k+1)$ is a $D_n$-module
with a basis given by the basis graphs of $\Z[x] S_n \ek w_\j$.
In particular, we obtain
\begin{equation}\label{Dnmod}
I(k)/I(k+1)\ \cong\
\bigoplus_{\j\in P(n,k)} \Z[x] (S_n \ek w_\j +I(k+1))/I(k+1)
\end{equation}
As multiplication from the right by $w_\j$ commutes with the $D_n$-action,
it follows that each summand on the right hand side is isomorphic
to the module
\begin{equation}\label{Vnkdef}
\Vnk=\Z[x] S_n\ek+I(k+1)/I(k+1).
\end{equation}
Combinatorially, it is spanned by graphs with exactly $k$ horizontal
edges in the lower part, where the $i$-th edge connects the lower vertices
$2i-1$ and $2i$.
As additional notation, let $s_i=(i,i+1)$ be the transposition of
the numbers $i$ and $i+1$, and let $\WBk$ be the subgroup of $S_n$
generated by the elements $s_{2i-1}$, $1\leq i\leq k$ and
by $s_{2i}s_{2i-1}s_{2i+1}s_{2i}=(i,i+2)(i+1,i+3)$, $1\leq i<k$.
It is well-known
that $W(B_k)$ is isomorphic to the semidirect product of
$(\Z_2)^k$ with $S_k$. We have the following simple properties.

\begin{lemma}\label{basicVnk} (a) The $\Z[x]$-rank of $\Vnk$ is equal
to $n!/2^kk!$. Moreover, as an $S_n$-module, $\Vnk\cong \Z[x] (S_n/\WBk)$.

(b) $\Z[x] S_n\ek S_n$ is isomorphic to
$\Z[x] S_n\ek\otimes_{\Z[x]S_{2k+1,n}}\ek S_n$
as $\Z[x]S_n-\Z[x]S_n$-bimodule.

(c) The commutant of $D_n$ on $\Vnk$ is given by $\Z[x] S_{2k+1,n}$.

(d) The algebra $D_n$ is faithfully represented on
$\bigoplus_{k=0}^{[n/2]}\Vnk$.
\end{lemma}

$Proof.$ The first statement is straightforward to prove. The second
statement follows from the fact that $S_{2k+1,n}$, which leaves
the numbers 1 until $2k$ fixed, commutes with $\ek$, from which one
can deduce that $\Z[x] S_n\ek$ is a free  $\Z[x] S_{2k+1,n}$ right
module, and that  $\Z[x] \ek S_n$ is a free $\Z[x] S_{2k+1,n}$ left-module.
As to the
statement (c), it is easy to see that $\Z[x] S_{2k+1,n}$ is contained in
the commutant. As $\ek$ is
a cyclic vector for $\Vnk$, any element $b$ in the commutant of $D_n$
is already completely determined by its action on $\ek$.
It is easy to inspect by multiplying graphs that $\ek d$ is in
$\ek S_{2k+1,n}+I(k+1)$ for any $d\in S_n\ek$. Hence it follows
$$x^kb\ek = b\ek^2=\ek(b\ek)=\pi\ek$$
for some $\pi\in \Z[x] S_{2k+1,n}$.
To prove the last statement, we use the fact that
$\Q(x)\otimes _{\Z[x]}D_n$ is semisimple (see e.g. \cite{HW}).
Hence its left regular representation is faithful. But by the
discussion in this section, see \ref{Dnmod} and \ref{Vnkdef},
$D_n$ has a filtration of $D_n$-modules, each of whose factors
is isomorphic to a $\Vnk$. By semisimplicity, we can replace
this by a direct sum of modules each of which is isomorphic to a $\Vnk$.

\subsection{Decomposition}\label{decomposition}
In the following we are primarily interested in the $S_n$-action
on $\Vnk$. For simplicity, we do this over the ring $\Z$; the
results are exactly the same for the ring $\Z[x]$.
We shall need the decomposition of
the module $\Vnk$ as a $\Z S_{3,n}$-module, where $S_{3,n}$
is the group of permutations of letters 3 until $n$. In view
of the last lemma, it is clear that we obtain a decomposition
in terms of  $ S_{3,n}$-orbits of $S_n/\WBk$, i.e. in terms
of cosets $ S_{3,n}w\WBk$. We shall describe these double
cosets in terms of specially chosen elements $w$ whose
meaning will become clear later. If $i\leq j$, we shall use the
notation $s_{i,j}=s_is_{i+1}\ ...\ s_j$.
Not surprisingly,
the size of such double cosets depends on the intersection
$w^{-1}\{ 1,2\}\cap [2k+1,n]$. We list  the decomposition
of $\Vnk$ into $ S_{3,n}$-modules in the table below.

\vskip .4cm

\begin{tabular}{c|c|c|c}
 $w$
 & $\Z S_{3,n}w\WBk\cong$ & dim& $\#$ of modules\\
& & & \\
\hline
&  & &\\
$s_{2,j_2}s_{1,j_1}, j_1\geq 2k, j_2>2k$ & $V_{n-2}^{(k)}$ & $\frac{(n-2)!}{2^kk!}$&
$(n-2k)(n-2k-1)$  \\
&  & &\\
\hline
&  &  &\\
$s_{1,j_1}$ or $s_{2,j_2},\ j_1,j_2>2k$ & $V_{n-2}^{(k-1)}$ &$\frac{(n-2)!}{2^kk!}2k$& $ 2(n-2k)$ \\
&  & & \\
\hline
&  & & \\
1 & $V_{n-2}^{(k-1)}$ &$\frac{(n-2)!}{2^kk!}2k$& $ 1$ \\
&  & & \\
\hline
&  & & \\
(23)& $V_{n-2}^{(k-2)}$ &$\frac{(n-2)!}{2^kk!}2k(2k-2)$& $ 1$ \\
\end{tabular}

\subsection{Length function}\label{lengths}
Similarly as for elements in reflection
groups, one can define a length function for basis graphs of
the Brauer algebra. Recall that for a permutation $w\in S_n$, its length
$\ell(w)$ is the minimum number of factors in an expression of $w$
as a product of simple reflections; interpreting
$w$ as a graph as above, $\ell(w)$ would be the number of crossings in
that graph with the following $caveat$: The element $\ek$ is drawn
fixed and must be left unchanged; e.g. the element $s_1s_2es_2s_1$ has
length 4, even though the corresponding graph in the Brauer algebra could
be drawn without any crossings.
The precise definition of the length $\ell(d)$ of a basis graph
$d\in D_n$ with exactly $2k$
horizontal edges is given by
$$\ell(d)={\rm min}\{ \ell(w_1)+\ell(w_2),\ w_1\ek w_2=d,\ w_1,w_2\in S_n\}.$$
We will also call graphs of the form $w\ek$ basis graphs of
the module $\Vnk$.
For given $d$, there can be more than one $w$ with
$w\ek=d$ and $\ell(w)=\ell(d)$, e.g. $s_1s_2\ek=s_3s_2\ek$ for $k\geq 2$.
To pin down a specific choice,
it will be convenient to use the notation $s_{i,j}=s_is_{i+1}\ ...\ s_j$ for
$i<j$. It is well-known that the elements $w$ of $S_n$ can be written uniquely
in the form $w=t_{n-1}t_{n-2}\ ...\ t_1$, where $t_j=1$ or $t_j=s_{i_j,j}$
with $1\leq i_j\leq i$ and $1\leq i<n$. This can be easily seen as
follows: For given $w\in S_n$, there exists a unique $t_{n-1}$ such
that $t_{n-1}(n)=w(n)$. Hence $w'=t_{n-1}^{-1}W(n)=n$ and we can view
$w'$ as an element of $S_{n-1}$. The general claim now follows by
induction on $n$. We will apply a similar strategy for defining
basis elements for $\Vnk$.
Using the notation for the $t_j$'s,
we now define for $k\leq n/2$ the set
\begin{equation}\label{Bndef}
B_{n,k}\ =\ \{ (t_{n-1}t_{n-2}\ ... t_{2k}t_{2k-2}\ ...\ t_2\},
\end{equation}
Observe that $B_{n,k}$ has $n!/2^kk!$ elements.


\begin{lemma}\label{Vnklength} (a) The module $\Vnk$ has a basis
$\{ wv_1=v_{w\ek},\ w\in B_{n,k}\}$ with $\ell(w\ek)=\ell(w)$. Here $\ell(w)$
is the number of factors for $w$ in Def. \ref{Bndef}, and $v_1=\ek+I(k+1)\in
\Vnk$.

(b)
We have $|\ell(s_id)-\ell(d)|\leq 1$ for any basis graph for $\Vnk$.
Equality of lengths holds only if $s_id=d$.

(c) Let $S_3^{(i)}$ be the subgroup of $S_n$
generated by $s_i$ and $s_{i+1}$. Then each  $S_3^{(i)}$-orbit
$O$ of basis graphs in $\Vnk$ has the order structure of $S_3^{(i)}/W$,
where $W$ is either the trivial subgroup or the subgroup generated
by $s_i$ or by $s_{i+1}$.

\end{lemma}

$Proof.$ Let $d=w\ek$ be a graph in $S_n\ek$. Using exactly the
same arguments as given before Def \ref{Bndef}, we determine
$t_{n-1},\ ...,\ t_{2k}$ such that $d'=(t_{n-1}t_{n-2}\ ...\ t_{2k})^{-1}d$
is a graph in $S_{2k}\ek$., i.e. $d'$ can be viewed as
a graph in $D_{2k}$ with only horizontal edges to which we add
$n-2k$ strictly vertical edges to the right.
Let $i_{k-2}$ be the label of the upper vertex of $d'$ which is
connected with the upper $2k$-th vertex and set $t_{2k-2}=s_{i_{k-2},k-2}$.
Then the upper $2k$-th and $(2k-1)$-st vertices of $d''=t_{2k-2}^{-1}d''$
are connected by a horizontal edge. Proceeding in this way, we eventually
transform $d$ into $\ek$. Hence every graph in $S_n\ek$ can be written
as $w\ek$, with $w\in B_{n,k}$.

To show that the $w$ constructed in the last paragraph has minimal
length, let $v\in S_n$ be such that $v\ek = d$.
Let $1\leq r\neq s\leq k$. Then it is easy to see, e.g. by drawing
pictures, that we have at least zero, one or two intersections
between edges emanating from $2r-1, 2r, 2s-1, 2s$ for $1\leq r<s\leq k$ if
$[v(2r-1),v(2r)]\cap [v(2s-1),v(2s)]$ is empty,
is a proper subintervall of both intervals, or
is equal to one of the two intervals respectively.
Let us call this minimum number $c(r,s)$.
Moreover, we get an additional crossing for each inversion,
i.e. for each pair
$1\leq a<b\leq n$ with $b>2k$ for which  $v(a)>v(b)$.
It is not hard to check that the number of inversions
(with $b>2k$) is independent of the choice of $v$.
Hence
\begin{equation}\label{lengthformula}
\ell(d)\geq\sum_{1\leq a<b,\ b>2k}inv(a,b)+
\sum_{1\leq r<s\leq k}c(r,s),
\end{equation}
where
$inv(a,b)$ is equal to 0 if $v(a)<v(b)$ and equal to 1 if $v(a)>v(b)$.
It remains to check that  the right hand side is equal to
$\ell(t_{n-1}t_{n-2}\ ...\ t_{2k})+\ell(t_{2k-2}t_{2k-4}\ ...\ t_2)=
\ell(w)$ for $w$ as constructed in the previous paragraph.
This is easy. Hence we have equality in Eq. \ref{lengthformula}.

Part (b) can now be checked in a fairly straightforward way using
the explicit formula for the length.
Also part (c) is either known from the symmetric group case, or
it can be checked in a straightforward way. E.g. if the numbers
$i$, $i+1$ and $i+2$ label vertices
belonging to three different horizontal upper
edges of $d$, say $(i,j_1)$, $(i+1,j_2)$ and $(i+2,j_3)$, the
action of $S_3(i)$ results in permuting the second coordinates,
and it is easy to see that the lowest element is given
if $j_1<j_2<j_3$. In this case, it can be explicitly checked,
for instance via pictures, that the map $w\mapsto w(i,j_1)(i+1,j_2)(i+2,j_3)$
is order-preserving. The other cases are similar and easier.

\subsection{Braids and Hecke algebras}
Recall that Artin' s braid group $AB_n$ is defined via
generators $\sigma_i$, $1\leq i\leq n-1$ and relations
$\sigma_i\sigma_j=\sigma_j\sigma_i$ for $|i-j|>1$ and
$\sigma_i\sigma_{i+1}\sigma_i=\sigma_{i+1}\sigma_i\sigma_{i+1}$.
It will also be convenient to introduce the notation
$\sigma^+_{k,l}=\sigma_k\sigma_{k+1}\ ...\ \sigma_{l-1}\sigma_l$ if $k<l$ and
$\sigma^+_{k,l}=\sigma_k\sigma_{k-1}\ ...\ \sigma_{l+1}\sigma_l$ if $k>l$.
Similarly, the expressions $\sigma^-_{k,l}$ are defined as above,
with $\sigma_i$ replaced by $\sigma_i^{-1}$ for $k\leq i\leq l$.
Similarly, one defines elements $\gp_{k,l}$ and $\gm_{k,l}$
in terms of the generators $g_i$ of the Hecke algebra (see below).
We have the following simple lemma, which is easy to prove.

\begin{lemma} \label{braidrelations}
(a) The map $\Phi : \sigma_i\mapsto  \sigma_{2i}\sigma_{2i+1}
\sigma_{2i-1}^{-1} \sigma_{2i}^{-1}$
defines a homomorphism of the braid group $AB_n$ into $AB_{2n}$.

(b) $\sigma_j^{\pm 1}\sigma^{\pm}_{k,l}=
\sigma^{\pm}_{k,l}\sigma^{\pm 1}_{j-1}$ if $k<l$ and $k<j\leq l$.

(c)  $\sigma_j^{\pm 1}\sigma^{\pm}_{k,l}=
\sigma^{\pm}_{k,l}\sigma^{\pm 1}_{j+1}$ if $l<k$ and $l\leq j<k$.
\end{lemma}

The Hecke algebra $H_n$ of type $A_{n-1}$
is the $\Zq$-algebra defined by generators $g_i$, $1\leq i<n$ and relations
$g_ig_{i+1}g_i=g_{i+1}g_ig_{i+1}$ and $g_ig_j=g_jg_i$ for $|i-j|>1$.
It has a basis $(g_w)_{w\in S_n}$ such that
\begin{equation}
g_ig_w=
\begin{cases} g_{s_iw}  & \text{if $\ell(s_iw)>\ell(w)$,}\\
(q-1)g_w + qg_{s_iw}
&\text{if  $\ell(s_iw)<\ell(w)$. }
\end{cases}
\end{equation}

It will be convenient to define the module $\Vnk$
as a $\Zq$-module
with a basis $(v_d, d=w\ek, w\in B_{n,k})$.
We will subsequently define actions of the Hecke algebra
and of a $q$-deformation of the  Brauer algebra on this module
which will specialize to the known actions if we restrict
to the classical Brauer algebra. So no confusion should arise
from this slight abuse of notation.
We now define an action of the generators $g_i$ of $H_n$ on $\Vnk$
as follows:

\begin{equation}\label{Heckeaction}
g_iv_d=
\begin{cases} qv_d  & \text{if $s_id=d$,}\\
v_{s_id}  & \text{if $\ell(s_id)>\ell(d)$,}\\
(q-1)v_d + qv_{s_id}
&\text{if  $\ell(s_id)<\ell(d)$. }
\end{cases}
\end{equation}

\begin{proposition} The action defined in \ref{Heckeaction} makes
the $\Zq$-module $\Vnk$ into an $H_n$-module.
\end{proposition}

$Proof.$ This could be checked by identifying  $\Vnk$
with a quotient of $H_n$, see Lemma \ref{quotientmodule}.
Here we check the relations directly as follows:  For given $g_i$ and $g_{i+1}$,
this only needs to be done on the subspaces spanned by the
$S_3(i)$-orbits of the basis graphs. These are either 6 or 3-dimensional.
As the definition of the action only depends on the order structure
of the basis elements, it follows from Lemma \ref{Vnklength}
that the actions on these subspaces coincides with the left regular
representation of $H_3(i)$ in the 6-dimensional case, and with
a representation on a coset space in the 3-dimensional case.
It is not hard to check that in the latter case we obtain the same
matrices as the ones for $g_1$ and $g_2$ in Section \ref{lowdim}.
The relation $g_ig_j=g_jg_i$ can be checked in a similar way and is easier.

\medskip
Let $1\leq i\leq j$ and let
$n,m \geq j$. We will later need the following relations, which can be
proved by straightforward calculations, similar to the ones in
Lemma \ref{braidrelations}.
\begin{equation}\label{Heckerel1}
\gp_{i,n}\gm_{j,m}\ =\
\begin{cases} \gm_{j+1,m+1}\gp_{i,n}  & \text{if $m<n$,}\\
\gm_{3,n}\gp_{2,n-1}
&\text{if  $m\geq n$. }
\end{cases}
\end{equation}

\begin{equation}\label{Heckerel2}
\gp_{i,n}\gp_{j,m}\ =\
\begin{cases} \gp_{j+1,m+1}\gp_{i,n}
& \text{if $m<n$,}\\
(q-1)\gp_{j+1,n}\gp_{i,m} + q\gp_{i,m}\gp_{i,n-1}
&\text{if  $m\geq n$. }
\end{cases}
\end{equation}
Moreover, the same relations hold if we simultaneously replace
all $+$ signs with $-$ signs and vice versa, in each of the formulas
above.

\subsection{Other versions} Obviously, we also obtain other $S_n$-modules
in the Brauer algebra after conjugating $\ek$ by a permutation.
These modules can be generalized to Hecke algebra modules as before.
However, as already remarked at the beginning of Section \ref{lengths},
we may get different length functions for the resulting graphs.
We deal here with the special case where  $\ek$ is replaced by the
same graph except that the two leftmost horizontal edges are replaced
by two vertical edges to keep notation simpler.
We denote this element by $e_{(2,k)}$.
Similarly, we can also define the module $V^{(2,k)}_n$ both for
the Brauer algebra, and for the Hecke algebra; we denote
the vector corresponding to the element  $e_{(2,k)}$
by  $v_1^{(2,k)}$. The length function
for basis elements of the module $V^{(2,k)}_n$ is defined as
before for $\Vnk$, except that $\ek$ is replaced by  $e_{(2,k)}$.
We shall need the following technical lemma:

\begin{lemma}\label{quotientmodule}
Let $L=L_{(n,k)}$ be the left ideal in $H_n$
generated by $g_{2i-1}-q,\ 1\leq i\leq k$
and by $g_{2i+1}g_{2i}-g_{2i-1}g_{2i},\ 1\leq i\leq k-1$.
Then $H_n/L$ is a free $Z[q,q^{-1}]$-module of rank $n!/2^{k}k!$
\end{lemma}

$Proof.$  For $k=0$, the module
$V_n^{(0)}$ is just the left regular representation
of $H_n$, and there is nothing to show. If $k>0$, it follows
from the definitions that  $L$ is contained
in the annihilator of the vector $v_1$. Hence $H_n/L$ has at least
dimension $n!/2^kk!$. So
it suffices to show that $H_n$ is equal to the span of
$B_{n,k}$ and $L$.
We shall show this by induction on $n$ and $k$.
Let us first show that it suffices to prove this for $n=2k$.
Indeed, in this case the claim for  $n>k$ follows by induction on $n$
by observing that
\begin{equation}
H_{n+1}=\bigoplus_{a=1}^{n+1} g_{a,n}H_n=
span \bigcup_{a=1}^{n+1}g_{a,n}B_{n,k}\cup g_{a,n}L^{(r,k)}_n,
\end{equation}
where we set $g_{n+1,n}=1$.
The claim now follows from the fact that
$B_{n+1,k}=\bigcup_{a=1}^{n+1}g_{a,n}B^{n,k}$, see Def. \ref{Bndef}.

It remains to show the claim for $n=2k$ and $r=0$, which we again do
by induction on $k$, with $k=1$ being trivially true.
By the above, the claim also holds for $n=2k+1$, with
$B^{(k)}_{2k+1}=\bigcup g_{a,2k}B^{(k)}_{2k}$.
Let $b=g_{i_{2k},2k}b'\in B^{(k)}_{2k+1}$.
If $g_{i_{2k},2k}=1$ then we have
\begin{equation}\label{case1}
g_{i_{2k+1},2k+1}b-qg_{i_{2k+1},2k}b\in L
\end{equation}
while if  $g_{i_{2k},2k}\neq 1$, we have
\begin{equation}\label{case2}
g_{i_{2k+1},2k+1}b-g_{i_{2k+1},2k}g_{i_{2k-1},2k}b'g_{2k-1}g_{2k}
= g_{i_{2k+1},2k}g_{i_{2k-1},2k}b'(g_{2k+1}g_{2k}-g_{2k-1}g_{2k})\in L
\end{equation}
It follows that the elements in \ref{case1} and \ref{case2}
together with the ones in $L^{(2k)}_{2k+1}$ and the ones in
$B^{(k)}_{2k+1}= B^{(k+1)}_{2k+2}$ span $H_n$, as required.
\vskip .2cm
\begin{corollary}\label{quotientcor}
Let $L^{(r)}_{n,k}$ be the ideal generated by
the elements $g_{r+2i-1}-q$, $1\leq r\leq k$ and by
$g_{r+2i-1}g_{r+2i}-g_{r+2i+1}g_{r+2i}$, $1\leq r\leq k-1$.
Then again $H_n/L^{(r)}_{n,k}$ is a free $\Z[q,q^{-1}]$-module
of rank $n!/2^kk!$
\end{corollary}

$Proof.$
Conjugating the ideal $L_{n,k}$ by
the element $g_{1,2k}g_{2,2k+1}\ ...\ g_{r, 2k+r-1}$
gives us the ideal  $L^{(r)}_{n,k}$.

\subsection{$H_{3,n}$-modules}
We can now use these results to define certain $H_{3,n}$-module
morphisms in $\Vnk$ which will be needed later. First of all,
we replace the elements $w$ in the table of Section \ref{decomposition}
by elements $g_w$ by
replacing
$s_{2,j_2}$ by $g^+_{2,j_2}$ and replacing $s_{1,j_1}$ by
$g^-_{1,j_1}$. Then we can show the following:

\begin{lemma}\label{H3nmorphisms} Let $w$ be an element as in the table
of Section \ref{decomposition}, and let $g_w$ be as just defined.
Then we get a decomposition $\Vnk\cong \bigoplus H_{3,n}g_wv^{(k)}_1$
as $H_{3,n}$-modules analogous to the one in  Section \ref{decomposition}.
In particular, we have the following well-defined
$H_{3,n}$ homomorphisms:

(a) $hg^-_{1,j_1}v^{(k)}_1\mapsto hg^-_{3,j_1}v^{(k)}_1$, $h\in H_{3,n}$, $j_1>2k$,

(b) $hg^+_{2,j_2}v^{(k)}_1\mapsto hg^+_{3,j_2}v^{(k)}_1$, $h\in H_{3,n}$, $j_2>2k$,

(c) $hg_2v^{(k)}_1\mapsto hv^{(k)}_1$, $h\in H_{3,n}$,

(d) $hg^+_{2,j_2}v^{(k)}_1\mapsto hg^-_{2,j_2}v^{(k)}_1$ and
 $hg^+_{2,j_2}v^{(k)}_1\mapsto hg^-_{2,j_2}v^{(k)}_1$,  $h\in H_{3,n}$.
\end{lemma}

$Proof.$ The only nontrivial part in the proof is to show that the maps
are well-defined. Observe that
 in case (c) the annihilator of $g_2v^{(k)}_1$ in $H_{3,n}$
contains the elements
$g_{2i-1}-q$, $3\leq i\leq k$ and $g_{2i+1}g_{2i}-g_{2i-1}g_{2i}$,
$3\leq i< k$. By Lemma \ref{quotientmodule}, the quotient of
$H_{3,n}$ with the left ideal $L$ generated by these elements has
rank $(n-2)!/2^{k-2}(k-2)!$, which coincides with the rank of
the module $H_{3,n}g_2v^{(k)}_1$, see the table in Section
\ref{decomposition}. Hence the annihilator coincides with $L$,
which is obviously contained in the annihilator of $v_1^{(k)}$.
It follows that the homomorphism is well-defined.
One similarly determines annihilator ideals in the other cases,
using
Lemma \ref{quotientmodule}, Corollary \ref{quotientcor}
and the table in Section \ref{decomposition}. The claim follows
as before.

\ignore{To show (a) and (b), one first observes that
the dimensions of the $H_{3,n}$-modules generated by the various
cyclic vectors are given in the table in Section \ref{decomposition}.
It then is not hard to show that the annihilators of the cyclic
vectors on both sides of the maps in (a) and (b) coincide,
using Lemma \ref{quotientmodule} and Corollary \ref{quotientcor}.
Also, for (d) it is easy to calculate that the annihilators for
 $hg^\pm_{2,j_2}v^{(k)}_1$ and  $hg^\pm_{1,j_1}v^{(k)}_1$
are generated by $g_{2i}-q$, $2\leq i\leq k$ and by
$g_{2i}g_{2i+1}-g_{2i+2}g_{2i+1}$, by the same arguments as before.
The isomorphism follows from that.}

\section{Deformation of module tensor categories}\label{motivation}

\subsection{Motivation and deformation conditions} This and the subsequent
subsection only serve to motivate the following
definitions. They are less self-contained and less rigorous
than the other parts of this paper, which can be read independently
of this section. For background for categorical notions see
e.g. the book \cite{Kassel} and references therein, and the paper \cite{Os}.

It is well-known that for groups $H\subset G$, we can make the representations
of $H$ into a module category of $Rep(G)$.  The right module action is
defined for  $V$  an $H$-module, $W$  a $G$-module by
$V\otimes W = V\otimes Res(W)$, where $Res(W)$ is $W$
viewed as an $H$-module. In particular, we obtain
embeddings
\begin{equation}\label{tensormodule}
\End_H(V)\otimes \End_G(W)\to \End_H(V\otimes Res(W)).
\end{equation}

The idea for the construction of the new $q$-Brauer algebra can now
be stated very easily, which we will do on the level
of Lie algebras. Let $\h\subset\g$ be semisimple Lie algebras.
There exist canonical $q$-deformations of their universal enveloping
algebras due to Drinfeld and Jimbo. It is known that these deformations
usually are not compatible with the inclusion $\h\subset\g$.
Hence we weaken the problem and ask for a compatible deformation
of $Rep(\h)$ as a module category over $Rep(\g)$. More precisely,
we require the following conditions:

{\it (A) Same restriction rules:}
If $\Ca$ is the (finite-dimensional)  representation category of a
Drinfeld-Jimbo quantum group corresponding to $\g$,
we would like to find a module  category $\Da$
with the same Grothendieck semigroup as $Rep(\h)$ and with a
right tensor module
action as in \ref{tensormodule} which
should be compatible with the identifications of Grothendieck
semigroups.

\ignore{
 In particular, we get for any object $X$ in
$\Ca$ canonical embeddings
$$\End_\Ca(X)\cong \End_{\Da}(\1)\otimes \End_\Ca(X)\to
\End_\Da(\1\otimes X)\cong \End_\Da(X),$$
where $\1$ is the trivial object of $\Ca$, which is required
to remain simple if viewed as an object of $\Da$.}

{\it (B) Compatible traces}
In addition  $\Ca$ is a spherical category,
i.e. it has canonical duality morphisms which lead to canonical
traces for $\End(X)$, for any object $X$ in $\Ca$
(see e.g. the chapter on duality in \cite{Kassel}). We also require
that these extend in a compatible way to our module category.
This condition is equivalent to a fundamental notion in the study
of subfactors known as the commuting square condition.
We will state it in this context as follows:

In a spherical category, there exists for every object
$Z$ in $\Ca$ a canonical trace $Tr_Z$ on $\End_\Ca(Z)$;
we will denote by $tr_Z$ the multiple of $Tr_Z$ such that
$tr_Z(1)=1$.
We now require extensions of $Tr_Z$
to $\End_\Da(Z)$ such that the following holds:
\begin{equation}\label{commsquare}
E(a)\in \End_\Ca(X) \quad {\rm for \ any\ }a\in \End_\Ca(X\otimes Y),
\ X,Y\in Ob(\Ca);
\end{equation}
here $E$ is the orthogonal projection onto the subalgebra $\End_\Da(X)
\cong \End_\Da(X)\otimes 1\subset \End_\Da(X\otimes Y)$ with
respect to the bilinear form $(b,c)=tr(bc)$;
for more details see Section \ref{SecSemisimp}.

\subsection{Some relations}
We give some examples how Cond. \ref{commsquare} forces
relations for a deformation of Brauer's centralizer algebra,
if we take for $\g=\sl_N$ and for the subalgebra $\h=\so_N$, with $N$
odd to avoid needless complications. We denote by $V$ the object
corresponding to the vector representation of $\sl_N$ resp. of $\so_N$
both in $\Ca$ and in the module category $\Da$.
It is well-known that $\End_\Ca(V^{\otimes n})$ is generated by
a representation of the Hecke algebra $H_n$.  We shall denote the images
of the generators again just by $g_i$. The canonical traces mentioned
before are known under the name Markov traces; see Section \ref{SecMarkov}
for details. In this context, Cond. \ref{commsquare} translates
for $X=V^{\otimes n}$ and $Y=V$ to the condition
\begin{equation}\label{qBrMarkov}
tr(bg_n)=tr(b)tr(g_n),\quad b\in \End_\Da(V^{\otimes n})\otimes 1.
\end{equation}

Let $\eb$ denote the projection in $\End_\Da(X^{\otimes 2})$
onto the object in $X^{\otimes 2}$
corresponding to the trivial representation of $\so_N$, which is
a subrepresentation of the symmetrization of the vector representation.
One deduces from this that $\eb g_1=qg_1$, as the eigenprojection of
$g_1$ with eigenvalue $q$ projects onto the object corresponding to
the symmetrization of the vector representation.

We shall also denote the embedding $\eb\otimes 1$ of $\eb$ into
$\End_\Da(X^{\otimes 3})$ just by $\eb$. Then $\eb$ also projects
onto a simple object in $X^{\otimes 3}$, and hence $\eb g_2\eb=\al \eb$
for a scalar $\al$. To calculate this scalar, we use the requirements
concerning the conditional expectation: By definition, $E(\eb g_2)$
is the unique element in $\End_\Da(X^{\otimes 2})$ such that
$tr_{X^{\otimes 3}}(a\eb g_2)=tr_{X^{\otimes 2}}(aE(\eb g_2))$
for all $a\in \End_\Da(X^{\otimes 2})$.
It follows from Eq. \ref{qBrMarkov} and \ref{condexp}
that the solution is
$E(\eb g_2)=tr(g_2)\eb$.
But then we also have
$$tr_{X^{\otimes 3}}(\al \eb)=tr_{X^{\otimes 3}}(\eb g_2\eb)=
tr_{X^{\otimes 3}}(\eb \eb g_2)=tr_{X^{\otimes 2}}(tr(g_2)\eb).$$
 Hence $\al=tr(g_2)$. Choosing suitable normalizations,
it is not hard to derive from these arguments the additional relations
$(E1)$ and $(E2)$ of the definition in the next section,
with $tr(g_2)=q^N/[N]$ and
$e=[N]\eb$ (see next section for notations).
Moreover, we will check later that the
condition \ref{commsquare} holds if we also have relation $(E3)$.

\begin{remark}\label{forcerelation}
It is possible to derive relation $(E3)$ in Section \ref{qBrdef}
from
condition \ref{commsquare} and relations $(H), (E1)$ and $(E2)$.
More precisely, these conditions and relations
essentially determine the matrices
of $g_3$ in all irreducible representations of $\qB_4$ with respect
to the path basis, see e.g. \cite{WHe} (for the Hecke algebra part)
and \cite{RW}. From this one can check that relation $(E3)$ has to
be satisfied as well. The proof is not very instructive, so we do
not give the details here.
\end{remark}

\section{$q$-Brauer algebras}\label{qBraueralgebras}

\subsection{Definitions}\label{qBrdef}
 Fix $N\in \Z\backslash\{0\}$ and let $[N]=(1-q^N)/(1-q\in\Z[q,q^{-1}]$.
The $q$-Brauer algebra $\BnN$ is defined over the ring $\Z[q,q^{-1}]$
 via generators
$g_1,g_2,\ ...\ g_{n-1}$ and $e$ and relations
\begin{enumerate}
\item[(H)] The elements $g_1,g_2,\ ...\ g_{n-1}$ satisfy the relations
of the Hecke algebra $H_n$.
\item[(E1)] $e^2=[N]e$,
\item[(E2)] $eg_i=g_ie$ for $i>2$, $eg_1=qe$, $eg_2e=q^Ne$ and
$eg_2^{-1}e=q^{-1}e$.
\item[(E3)] $g_2g_3g_1^{-1}g_2^{-1}e_{(2)}=e_{(2)}=
e_{(2)}g_2g_3g_1^{-1}g_2^{-1}$,
where $e_{(2)}=e(g_2g_3g_1^{-1}g_2^{-1})e$.
\end{enumerate}
We shall need a second version of the $q$-Brauer algebra, denoted by
$\qB_n(r,q)$ or just $\qB_n$ by carrying the information of the parameter
$N$ in the variable $r=q^N$. More precisely, the algebra $\qB_n(r,q)$
is defined over the ring $\Ra=\Z[q^{\pm 1}, r^{\pm 1}, (r-1)/(q-1)]$
via the same generators as before, with relations $(H)$ and $(E3)$
unchanged, and with
\begin{enumerate}
\item[$(E1)'$] $e^2=\frac{r-1}{q-1}e$,
\item[$(E2)'$] $eg_i=g_ie$ for $i>2$, $eg_1=qe$, $eg_2e=re$ and
$eg_2^{-1}e=q^{-1}e$.
\end{enumerate}

\begin{remark} 1. It should be clear that we get back the algebra
$\BnN$ from $\qB_n(r,q)$ by setting $r=q^N$.
In particular, we can use this to also define $\qB_n(0)$ as one of those
specializations, where
the direct definition would cause some (presumably minor) technical 
difficulties (see e.g. Lemma \ref{Brrelations},(g));  the author would like 
to thank Dung Tien Nguyen for pointing this out to him.
It is also easy to see that we get
the Brauer algebra $D_n(N)$ for $r=q^N$ in the limit $q\to 1$.
In this case $g_i$ becomes the simple reflection $s_i$ and
the element $e$ can be identified with the graph $e_{(1)}$.
In general, we prefer the algebra $\BnN$ as its defining ring is
more natural, and it is closer to the intended applications.
However, as the algebra $\qB_n(r,q)$ is generically semisimple,
it is sometimes more convenient to work with.
In many cases, the proofs are the same for both versions and we
will only give them for one version, sometimes without explicitly
mentioning the other version.

2. It is easy to see that the assignment
$g_i\mapsto g_i^T=g_i$ and $e\mapsto e^T=e$
defines a linear anti-automorphism $a\mapsto a^T$ of
$\qB_n(r,q)$. Similarly, the map $g_i\mapsto g_i^*=g_i^{-1}$
and $e\mapsto e^*=q^{1-N}e$ defines an anti-linear antiautomorphism
with respect to the involution of the ring $\Ra$  defined by $\bar q=q^{-1}$
and $\bar r=r^{-1}$.

3. We shall later show that the subalgebra of $\BnN$ resp. $\qB_n(r,q)$
generated
by the generators $g_1,g_2,\ ...\ g_{n-1}$ is indeed isomorphic to $H_n$.
If the reader feels uncomfortable with this, he should use different
notation for the generators of the Hecke algebras.

4. It may be instructive to some readers to visualize the relations via
graphical calculus for ribbon tensor categories (see e.g. \cite{Kassel},
\cite{Turaev}), with $e$ given by the composition $\cup\circ\cap$,
and $g_i$ given by a standard braid generator $\sigma_i$.
While this may give a somewhat better intuitive feel about the relations,
it does not provide a topological interpretation for our algebra.
E.g. in this usual tangle interpretation, $e_{(2)}$ would describe
the same topological object as $g_2g_3g_1g_2e_{(2)}$, while it can
be checked that these are different elements in $\qB_4$. It would
be interesting if one could find a topological interpretation of
our algebra.
\end{remark}

\subsection{Low-dimensional examples}\label{lowdim}
One checks directly for $n=2$ that $\BnN$ is spanned by the elements
$1, g_1$ and $e$. If $n=3$ one also easily shows that
$\BnN$ is spanned by the basis elements
$g_w$ of $H_3$ and the elements $h_1eh_2$, where $h_1\in\{ 1, g_2, g_1g_2\}$
and $h_2\in\{1, g_2, g_2g_1\}$. Hence its rank is at most 15.
On the other hand, consider the assignments
\begin{equation}\label{matrices1}
 g_1\mapsto
\begin{pmatrix}q&0&0\\
               0&0&q\\
               0&1&q-1
\end{pmatrix},
\quad  g_2\ \mapsto\
\begin{pmatrix}0& q&0\\
               1&q-1&0\\
               0&0&q
\end{pmatrix}
\quad {\rm and }\quad
e_1\mapsto
\begin{pmatrix}\frac{r-1}{q-1}& r&rq\\
               0&0&0\\
               0&0&0
\end{pmatrix}.
\end{equation}
It is easy to check that these matrices define a representation of
$\qB_3(r,q)$ whose image is a free $\Ra$-module of rank 9.
By calculating the determinant of
the matrix formed from the nonzero rows of
the matrices representing the elements $e$, $eg_2$ and $eg_2g_1$,
one can also determine for which algebraic relations for $r$ and $q$
this representation is not semisimple.
 We have the following Lemma:

\begin{lemma}\label{representations}
(a) The algebra $\qB_3(r,q)$ is a free $\Ra$-module of rank 15.

(b) We obtain a representation of $\qB_4(r,q)$ from the representation
in (a) by assigning to $g_3$ the matrix of $g_1$. It is equivalent to
the representation of $\qB_4(r,q)$ on $\qB_4(r,q)e(g_2g_3g_1^{-1}g_2^{-1})e$.
In particular, the ideal generated by $e_{(2)}$ has rank 9.

(c) We also have $e_{(2)}=e(g_2g_3g_1^{-1}g_2^{-1})e=
e(g_2^{-1}g_3g_1^{-1}g_2)e=e(g_2^{-1}g_3^{-1}g_1g_2)e$.
\end{lemma}

$Proof.$ We have already shown part (a). The fact that we also obtain
a representation of $\qB_4(r,q)$ as described in (b) is almost immediate.
It only remains to show that  $\qB_4(r,q)e_{(2)}$ is spanned by
$e_{(2)}$, $g_2e_{(2)}$ and $g_1g_2e_{(2)}$, which follows from the
$\qB_3(r,q)$ case and relation $(E3)$. Part (c) can be shown
by a direct calculation using $(E2)$, $g_i^{-1}=q^{-1}g_i+(q^{-1}-1)$ and
$g_i=qg_i^{-1}+(q-1)$ as well as the identity
$q^{-1}(q-1)eg_2e+q(q^{-1}-1)eg_2^{-1}e+(r-1)(q^{-1}-1)e=0$.

\subsection{Elements $\ek$}
In the following, we define elements $\ek$ in $\BnN$
inductively by $e_{(1)}=e$ and by
\begin{equation}\label{ekdefinition}
e_{(k+1)} \ =\  e\gp_{2,2k+1}\gm_{1,2k}\ek\ =\ e\Phi(\gp_{1,k})\ek\ =\
\ek\Phi(\gp_{k,1})e
\end{equation}
where $\Phi$ is defined as in Lemma \ref{braidrelations}
with $\sigma_i$s replaced by $g_i$s.
The equivalence of these and additional
expressions for $\ek$ will be proved in the following lemma.
For $q=1$, it is not hard to show that both definitions
produce the same graph in the usual Brauer algebra. The following
lemmas will indicate how the Brauer relations will extend to
these new algebras.

\begin{lemma}\label{Brrelations}  (a) The elements $\ek$ are well-defined.

\ignore{ (a) There exist $e_{(j,k)},e_{(k,j)}\in\BnN$ such that
$e_{(j)}e_{(j,k)}
=\ek=e_{(k,j)}e_{(j)}$ for $j<k$.}

(b) $g_{2j-1}g_{2j}\ek=g_{2j+1}g_{2j}\ek$ and
$g_{2j-1}^{-1}g_{2j}^{-1}\ek=g_{2j+1}^{-1}g_{2j}^{-1}\ek$
for $1\leq j < k$.

(c) $\gp_{1,2l}\ek=\gp_{2l+1,2}\ek$ and  $\gm_{1,2l}\ek=\gm_{2l+1,2}\ek$
for $l<k$,

(d) For any $j\leq k$ we have $e_{(j)}\ek=\ek e_{(j)}=[N]^j\ek$.

(e) $[N]^{j-1}e_{(k+1)}=e_{(j)}\gp_{2j,2k+1}\gm_{2j-1,2k}\ek$ for $1\leq j<k$.

(f) $e_{(j)}g_{2j}\ek=q^N[N]^{1-j}\ek$ for $1\leq j\leq k$.

(g) (\cite{Ng})  $\ek^T=\ek$ for $N\neq 0$ and $k\geq 1$.
\end{lemma}

$Proof.$ Part (a) is shown by induction on $k$, using the fact that
$\Phi(g_i)$ commutes with $e$ for $i>1$.
For part (b), the claim follows for $j=1$ from the definitions.
If $j>1$, we use $g_l\gp_{2,2k+1}\gm_{1,2k}=\gp_{2,2k+1}\gm_{1,2k}g_{l-2}$
for $l=2j-1,2j$, by Lemma \ref{braidrelations}, and induction assumption
to show the claim. Part (c) follows easily from (b) by induction on $l$.
For part (d), we use induction on $j$ and part (c) as
follows:
$$e_{(j+1)}\ek=[N]^{j-1}e\gp_{2,2j+1}\gm_{1,2j}\ek
= [N]^{j-1}e\gp_{2,2j-1}\gm_{2j-1,2}\ek = [N]^j\ek.
$$
\ignore{To get the last expression in the definition of $\ek$ it suffices
to prove the following formula by induction on $j$ for $j=1,2,\ ...,\ k-1$:
\begin{equation}\label{ejrecursion}
[N]e_{(j)}\gp_{2j,2k+1}\gm_{2j-1,2k}\ek=
e_{(j+1)}\gp_{2j+2,2k+1}\gm_{2j+1,2k}\ek.
\end{equation}
}
Part (e) is shown by induction on $j$ with $j=1$ being the
first definition of $\ek$. Moreover, we have
\begin{align} \notag
e_{(j+1)} \gp_{2j,2k+1}\gm_{2j-1,2k}\ek\ &=\
e\gp_{2,2j+1}\gm_{1,2j}\gp_{2j, 2k+1}\gm_{2j-1,2k}e_{(j)}\ek
\notag \\
& =\
[N]^j e\gp_{2,2k-1}\gm_{1,2k-2}\gp_{2j,2k+1}\gm_{2j-1,2k}
\ek = [N]^{j}e_{(k+1)}, \notag
\end{align}
which proves (e) using the induction assumption and part (d).
For part (f), observe that the left hand side of the statement is
equal to
$$eg^+_{2,2j}g^-_{1,2j-2}e_{(j-1)}\ek=
[N]^{j-1}eg^+_{2,2j}g^-_{2j-1,2}\ek
=[N]^{j-1}eg^-_{2j,3}g^+_{2,2j}\ek = q^N[N]^{j-1}\ek,$$
where we used (c), Lemma \ref{braidrelations} (b) and (c), and
relations (E2). Part (g) follows from the definitions and
Lemma \ref{representations},(c) for $k=1,2$, and by induction and part (e)
(with $j=k$) for $k>2$. The difficulty for $N=0$ and a complete proof in the
other cases was pointed out to the author by D. T. Nguyen in \cite{Ng}.

\begin{lemma}\label{reppreparation}
We have $e_{(j)}H_{n}\ek \subset H_{2j+1,n}\ek + \sum_{m\geq k+1}H_ne_{(m)}H_n$,
where $H_{r,s}$ is generated by $g_r, g_{r+1},\ ...,\ g_{s-1}$ and $j\leq k$.
Moreover, if
$j_1\geq 2k$ and $j_2\geq 2k+1$, we also have:

(a) $e\gp_{2,j_2}\gm_{1,j_1}\ek
=e_{(k+1)}\gp_{2k+2,j_2}\gm_{2k+1,j_1}$, if $j_1\geq 2k$ and $j_2\geq 2k+1$,

(b) $e\gp_{2,j_2}\gp_{1,j_1}\ek$ is equal to
$$e_{(k+1)}\gp_{2k+2,j_1}\gp_{2k+1,j_2}+
q^{N+1}(q-1)\sum_{l=1}^kq^{2l-2}(g_{2l+1}+1)\gp_{2l+2,j_2}\gp_{2l+1,j_1}\ek.$$ \end{lemma}

\ignore{
\begin{cases} q^N\gp_{j_1+3,j_2}\ek  & \text{if $j_1\leq 2k-2$ even
and $j_1<j_2-1$,}\\
\ek  & \text{if $j_1\leq 2k-2$ even
and $j_1=j_2-1$,}\\
q^{-1}\ek  & \text{if $j_1\leq 2k$ even
and $j_1=j_2$,}\\
q^{-2}\gm_{j_2+2,j_1}\ek  & \text{if $j_2\leq 2k-1$ odd
and $j_1>j_2$,}\\
e_{(k+1)}\gp_{2k+2,j_2}\gm_{2k+1,j_1}  & \text{if $j_1\geq 2k$ and $j_2\geq 2k+1$.}
\end{cases}
$$}

$Proof.$
We will use the analogous decomposition
of $H_n\ek$ into $H_{3,n}$-modules as in Section \ref{decomposition},
with the adjustments for the Hecke algebra case as explained at
the beginning of the next subsection. Let us first prove the claims for $j=1$.
Claim (a) follows almost immediately from  Lemma \ref{Brrelations}, (e).
This proves the first statement of the Lemma for elements
in the $H_{3,n}$ submodules in the first case of the table in Section
\ref{decomposition}. For submodules in the second case, the claim
follows from relations (E2), and the remaining cases are easy.

To prove part (a) for $j>1$
observe that the left hand side of the statement is
equal to
$$eg^+_{2,2j}g^-_{1,2j-2}e_{(j-1)}\ek=
[N]^{j-1}eg^+_{2,2j}g^-_{2j-1,2}\ek
=[N]^{j-1}eg^-_{2j,3}g^+_{2,2j}\ek = q^N[N]^{j-1}\ek,$$
where we used Lemma \ref{Brrelations}(c),
Lemma \ref{braidrelations} (b) and (c), and
relations (E2).
The first statement of the lemma for $j>1$
can now be done by a fairly straightforward induction on $j$,
using (a), Lemma  \ref{Brrelations}(e),(f) and the inductive definition
of $e_{(j)}$.

To prove part (b), we first observe that
\begin{equation}\label{brrel1}
e_{(j)}g_{2j}g_{2j+1}g_{2j-1}g_{2j}e_{(j)}=[N]^{j-1}
(q^2e_{(j+1)}+q^{N+1}(q-1)(g_{2j+1}+1)e_{(j)},
\end{equation}
which follows from $g_i=qg_i^{-1}+(q-1)$ and the relations proved so far.
We deduce from this
\begin{align}\label{brrel2}
e_{(j)}\gp_{2j,j_2}\gp_{2j-1,j_1}\ek\ =\ &
\frac{1}{[N]}q^2\gp_{2j+2,j_2}\gp_{2j+1,j_1}\ek +\\
 & + [N]^{j-1}q^{N+1}(q-1)(g_{2j+1}+1)\gp_{2j+2,j_2}\gp_{2j+1,j_1}\ek,
\notag
\end{align}
where we use  Lemma \ref{Brrelations}(e). This shows,
among other things that the first term
on the right hand side is an element in $\BnN$. We can now show
by downwards induction on $j$, starting with $j=k$, that

\begin{align}\label{brrel3}
\frac{1}{[N]^{j-1}}e_{(j)}\gp_{2j,j_1}\gp_{2j-1,j_1}\ek\ =\ &
\ \frac{1}{[N]^{j-1}}e_{(k+1)}\gp_{2k+2,j_1}\gp_{2k+1,j_2}\ek\ +\ \\
&\ q^{N+1}(q-1)\sum_{l=j}^kq^{2l-2}(g_{2l+1}+1)
\gp_{2l+2,j_2}\gp_{2l+1,j_1}\ek\;
\notag
\end{align}
This follows for $j=k$ almost immediately from Eq \ref{brrel2},
and for $j<k$ again from Eq \ref{brrel2} and induction assumption.
The desired identity now follows for $j=1$. We note again that even
though some of the expressions do not look like it, all the elements
involved are indeed in $\BnN$.

\medskip
\noindent
In analogy to the Brauer case, we can now define
\begin{equation}\label{I(j)def}
I(j)=\sum_{k=j}^{n/2} H_n\ek H_n.
\end{equation}
It follows from
the Lemma that $I(j)$ forms a two-sided
ideal in $\BnN$ for
$j=1,2,\ ...\ $ and we have the inclusions of two-sided ideals
$\BnN\supset I(1)\supset I(2)\supset\ ...$.

\begin{proposition}\label{dimestimate} The algebra $Br_n(r,q)$ is spanned
by $\sum_{k=0}^{n/2} H_n\ek H_n$. In particular, its dimension is
at most the one of the Brauer algebra $D_n$.
\end{proposition}

$Proof.$ To prove the first statement, it suffices to show that
the right hand side is invariant under multiplication by the
generators of $Br_n(r,q)$.
This is obvious for the Hecke algebra generators $g_i$.
It follows for left multiplication by $e$ from Lemma \ref{reppreparation}
for $j=1$. The same proof works for right multiplication, using
the involution $^T$, see Remark 2 after the definitions.

To prove the estimate for the dimension, observe that the annihilator
of $\ek$ in $H_n$, acting via left multiplication, contains
the left ideal  $L_{n,k}$ (see Lemma \ref{quotientmodule}).
Hence the dimension of $H_n\ek$ is at most equal to the
dimension of $\Vnk$, which is equal to the number of graphs
$S_n\ek$ in the Brauer algebra. One similarly shows
that the dimension of $\ek H_n$ is $\leq $ the number of graphs
in $\ek S_n$. Finally, it follows as in Lemma \ref{basicVnk}
that $H_n\ek H_n$ is a quotient of $H_n\ek \otimes_{H_{2k+1,n}} \ek H_n$,
where the latter has dimension $\leq \dim \Z[x] S_n\ek S_n$.
Hence the dimension of  $\qB_n(r,q)\cong \bigoplus H_n\ek H_n$
is at most the one of the Brauer algebra. This proves the other
inequality.

\subsection{The $\BnN$-module $\Vnk$} The results in the last section
show that $H_n\ek$ is a $\BnN$-module modulo $I(k+1)$.
We will show that it is isomorphic to
the Hecke algebra module $\Vnk$ after making it into a $\BnN$-module
by defining an action of $e$ on it. We will again use the decomposition of
$\Vnk$ into a direct sum of $H_{3,n}$-modules as in Lemma
\ref{H3nmorphisms} using the table in Section
\ref{decomposition}. As before, we will replace the elements
$s_{2,j_2}$ and $s_{1,j_1}$ in the first column of the table in
that section by $\gp_{2,j_2}$ and $\gm_{1,j_1}$ respectively
to obtain elements $g_w$ as before Lemma \ref{H3nmorphisms},
and we write
$\Vnk=\oplus_w H_{3,n}g_wv_1$ as a direct sum of
$H_{3,n}$-modules.
We now define the action of $e$ on $\Vnk$ by

\begin{equation}\label{eaction}
e h\gp_{2,j_2}\gm_{1,j_1}v_1\ =\
\begin{cases} q^Nh\gp_{3,j_2}v_1  & \text{if $\gm_{1,j_1}=1$,}\\
 q^{-1}h\gm_{j_2+1,j_1}v_1
& \text{if $\gp_{2,j_2}=1$,}\\
0 & \text{if $j_1\geq 2k$ and  $j_2\geq 2k+1$;}
\end{cases}
\end{equation}
moreover, we define  $ehv_1=[N]hv_1$ for $h\in H_{3,n}$
and $ehg_2v_1=q^nhv_1$ for the remaining two cases.
It follows from Lemma \ref{H3nmorphisms}
that the action of $e$ commutes with
the action of $H_{3,n}$; this implies that it is well-defined.
Moreover, observe that the image of $e$ on $\Vnk$ is equal to
$H_{3,n}v^{(k)}_1$. From this it follows easily that $eg_2e$ and $q^Ne$
act via the same map on  $\Vnk$; the same goes for
$eg_2^{-1}e$ and $q^{-1}e$. We have proved the following proposition,
except for part of relation $(E1)$ and $(E3)$, which will be proved
in the following subsections.

\begin{proposition}\label{representation}
The action of the elements $g_i$, $1\leq i<n$ and $e$ on $\Vnk$, as given in
Eq \ref{eaction} defines a representation of $\BnN$.
\end{proposition}

\subsection{Checking the relations $eg_1=qe=eg_1$}
As $e\Vnk=H_{3,n}v_1$, we see easily that the relation
$g_1e=qe$ is preserved. To check the relation $eg_1=qe$, we  express
the action of $e$ with respect to the original basis of the Hecke algebra
module $\Vnk$, which is now easier to do. Indeed, as we have already
established that $e$ commutes with $H_{3,n}$, it suffices to
calculate the action of $e$ on vectors of the form
$g^+_{2,j_2}g^+_{1,j_1}v_1^{(k)}$. It follows that
\begin{equation}\label{eact1}
e\gp_{2,j_2}\gp_{1,j_1}v_1^{(k)}=
q^{N+1}(q-1)\sum_{l=1}^k
q^{2l-2}(g_{2l+1}+1)\gp_{2l+2,j_2}\gp_{2l+1,j_1}v_1^{(k)}.
\end{equation}
This result holds for all $j_1\geq 2k$ and $j_2>2k$.
Moreover, observe that
\begin{equation}\label{eact2}
s_{2l+1}(s_{2l+2,j_2}s_{2l+1,j_1})>s_{2l+2,j_2}s_{2l+1,j_1}
\quad \Leftrightarrow \quad j_2\geq j_1,
\end{equation}
which holds for any $l\geq 0$ for which $2l+1 \leq {\rm min}(j_1,j_2)$.
We leave it to the reader to check, both for
$j_2\geq j_1$ and for $j_2<j_1$, using Eq \ref{eact1} and \ref{eact2} that
$$eg_1(\gp_{2,j_2}\gp_{1,j_1}v_1^{(k)})=
q^{N+1}(q-1)\sum_{l=1}^k
q^{2l-2}(g_{2l+1}+1)g_{2l+1}\gp_{2l+2,j_2}\gp_{2l+1,j_1}v_1^{(k)}.$$
The desired equality now follows from $(g_{2l+1}+1)g_{2l+1}=q(g_{2l+1}+1)$.

\subsection{Checking the relation $(E3)$}
Observe that $e_{(2)}\Vnk=H_{5,n}v^{(k)}_1$
by Lemma \ref{reppreparation}, from which one
easily deduces the first equation of relation (E3).
The second equality is more difficult to check.
We will first verify it for  $Br_4(r,q)$.
We then show that an arbitrary $\Vnk$ can be written
as a direct sum of  $Br_4(r,q)$-modules for each of
which relation $(E3)$ holds. This is done in several steps:

$Step\ 1:$ We show that relation $(E3)$ is preserved for $n=4$.
This is easy, as $e_{(2)}$ acts as 0 on $V_4^{(0)}$ and  $V_4^{(1)}$;
moreover, on  $V_4^{(2)}$, $g_1$ and $g_3$ act via the same linear
map, which also trivially implies that relation $(E3)$ is preserved.
It follows that $Br_4(r,q)$ has the same dimension as the Brauer
algebra $D_4$.

$Step\ 2:$ Let $\Brf$ be the algebra defined as $Br_4(r,q)$, except
for the relation $e_{(2)}=$ $e_{(2)}g_2g_3g_1^{-1}g_2^{-1}$. Observe
that we also have $e_{(2)}g_1=qe_{(2)}=e_{(2)}g_3$ in  $\Brf$.
As the subgroup generated by $s_1$ and $s_3$ has index 6 in $S_4$,
one deduces that  $e_{(2)}H_4$ is spanned by the elements
$e_{(2)}$, $e_{(2)}g_2$,
$e_{(2)}g_2g_1^{-1}$, $e_{(2)}g_2g_3, e_{(2)}g_2g_1^{-1}g_3$ and
$e_{(2)}g_2g_1^{-1}g_3g_2^{-1}$ in  $\Brf$. As $H_4e_{(2)}$ is
three-dimensional also in $\Brf$, it follows that $H_4e_{(2)}H_4$
has at most dimension 18 in $\Brf$. Now one checks directly for
the six spanning elements of  $e_{(2)}H_4$ that also in
$\Brf$ we have $e_{(2)}H_4e$ is spanned by $e_{(2)}$; indeed,
e.g. we have $e_{(2)}g_2g_1^{-1}g_3g_2^{-1}e=
eg_2g_1^{-1}g_3g_2^{-1}e_{(2)}$ (by definition of $e_{(2)}$),
which is equal to $e_{(2)}$ also in  $\Brf$. It follows from
this that also in $\Brf$ the ideal generated by $e_{(2)}$
is equal to $H_4e_{(2)}H_4$.

$Step\ 3:$ Let $v\in \Vnk$ and let $W=\Brf v$ be the $\Brf$-submodule
generated by $v$. Then $W$ is also a $Br_4(r,q)$-module if
$e_{(2)}g_2g_1^{-1}g_3g_2^{-1}v=e_{(2)}v$, $e_{(2)}g_2g_1^{-1}g_3v
= e_{(2)}g_2v$ and $e_{(2)}g_2g_3v=e_{(2)}g_2g_1v$. Indeed, if $I$
is the two-sided ideal of $\Brf$ such that $\Brf/I=Br_4(r,q)$, it is
easy to check that $I$ is generated by $e_{(2)}(g_2g_1^{-1}g_3g_2^{-1}-1)$,
$e_{(2)}g_2(g_1^{-1}g_3-1)$ and by $e_{(2)}g_2(g_3-g_1)$ as a
$\Brf$-left ideal. The claim follows from this and our assumptions.

\begin{lemma} The action of the generators of $Br_n$ on $\Vnk$ also
preserve relation $(E3)$.
\end{lemma}

$Proof.$ We decompose $\Vnk$ into a sum of cyclic $\Brf$-modules
of the form $W=\Brf v$, for suitable vectors $v$. It then suffices
to check that $W$ is also  a $Br_4(r,q)$-module by the criterion of
$Step\ 3$. Indeed, in this case $W\cong Br_4(r,q)/{\rm Ann}(v)$
is a $Br_4(r,q)$-module on which obviously also relation $(E3)$ holds.

The explicit checking of the criterion in $Step\ 3$ is somewhat tedious
as there are several different cases. It is easier to study the
combinatorics in the context of the original Brauer algebra.
Obviously, there is only something
to prove if $e_{(2)}H_4v\neq 0$. This implies that $k\geq 2$, and that
among the first four upper vertices at least
two belong to distinct horizontal edges. It remains to consider
the cases that 0, 1 or 2 of the first four upper vertices belong to
vertical edges.

Let us first consider the case with the first four upper vertices belonging
to four distinct horizontal edges. Multiplying such a graph by a suitable
permutation in $S_{5,n}$, if necessary, we can assume that each of these
four edges connect a vertex $\in \{1,2,3,4\}$ with a vertex
$\in \{ 5,6,7,8\}$ (observe that
 $S_{5,n}$ commutes with $D_4$, hence this multiplication induces
an isomorphism of $D_4$-modules). Among such graphs,
$s_4s_5s_6s_3s_4s_2e_{(4)}$ has the fewest crossings.
It will be convenient to pick the element
$v_0=g_4g_5g_6g_3^{-1}g_4^{-1}g_2^{-1}e_{(4)}$ in $Br_4(r,q)$.
We now leave it to the reader to check, using Lemma \ref{Brrelations}, (c)
that $e_{(2)}v_0=q^{N-1}v_0=e_{(2)}g_2g_3g_1^{-1}g_2^{-1}v_0$,
$e_{(2)}g_2v_0=q^{N-1}g_6v_0=e_{(2)}g_2g_3g_1^{-1}v_0$
and $e_{(2)}g_2g_3v_0=e_{(2)}g_5g_6v_0=e_{(2)}g_2g_1v_0$.

The case with three of the first four vertices connected to
three different horizontal edges, and the remaining one connected
to one of the lower row is done similarly. Here we can take
$v_0= g_4g_5g_6g_3^{-1}g_4^{-1}g_2e_{(2)}$, with essentially the
same calculations as before.

Next let us consider the case where two of the first four vertices
belong to horizontal edges which connect them with vertices larger than 4,
and that the other two vertices are connected to vertices in the lower row.
Again, it suffices to consider the cyclic module generated by the
element $v_0=g_4g_3g_5^{-1}g_4^{-1}g_2e_{(2)}$. Using the relations,
one checks that $e_{(2)}v_0=0=e_{(2)}g_2^{-1}g_1^{-1}g_3g_2v_0$,
$e_{(2)}g_2^{-1}v_0=q^{N-1}e_{(2)}=e_{(2)}g_2^{-1}g_1^{-1}g_3v_0$
and $e_{(2)}g_2^{-1}g_1^{-1}v_0=q^{N-1}g_5^{-1}e_{(2)}=
e_{(2)}g_2^{-1}g_3^{-1}v_0$.

In the remaining cases, we have at least two of the first four vertices
connected by a horizontal edge. We leave it to the reader to check
that these cases can be reduced to submodules generated by $ev_0$,
with $v_0$ as in one of the previous cases. This finishes the proof
of the lemma.

\subsection{Dimension} We can now prove the main theorem of this
section. We define for each basis graph $d$ of the Brauer
algebra $D_n$ an element $g_d\in \qB_n(r,q)$ as follows:
If $d$ has $2k$ horizontal edges, fix a reduced expression $d=w_1\ek w_2$
(see Section \ref{lengths}) with $w_1,w_2\in S_n$. Then we define
$g_d=g_{w_1}\ek g_{w_2}$; as usual, we abuse notation by denoting by
$\ek$ both a certain graph, and an element in $\qB_n(r,q)$.

\begin{theorem}\label{structuretheorem}
(a) The algebra $\qB_n(r,q)$ is a free
$\Z[q^{\pm 1}, r^{\pm 1}, (r-1)/(q-1)]$-module of rank $n!!=1\cdot 3\cdot\ ...
\ (2n-1)$ with basis $(g_d)$ labeled by the basis graphs of
the Brauer algebra.

(b) The algebra $\BnN$ is a free $\Z[q,q^{-1}]$-module of 
rank $n!!=1\cdot 3\cdot\ ... \ (2n-1)$ with spanning set
$(g_d)$ labeled by the basis graphs of the Brauer algebra.

(c) The algebra
$\qB_n(r,q)$ has the same decomposition into a direct
sum of simple matrix rings as a $\Q(r,q)$
algebra as the generic Brauer algebra $D_n$; this
also includes the restriction rules from, say, $\qB_n(r,q)$ to
$\qB_{n-1}(r,q)$, see Remark \ref{decompremark}.
\end{theorem}

$Proof.$  We have seen that there is
a faithful representation of the Brauer algebra $D_n$ on
$\bigoplus_{0\leq k\leq n/2}\Vnk$ in Lemma \ref{basicVnk}.
As this is a specialization of the representation of $Br_n(r,q)$
on the same direct sum of modules $\Vnk$, the dimension of
$Br_n(r,q)$ must be at least the one of $D_n$.

To prove the other inequality, observe that the annihilator
of $\ek$ in $H_n$, acting via left multiplication, contains
the left ideal  $L_{n,k}$ (see Lemma \ref{quotientmodule}).
Hence the dimension of $H_n\ek$ is at most equal to the
dimension of $\Vnk$, which is equal to the number of graphs
$S_n\ek$ in the Brauer algebra. One similarly shows
that the dimension of $\ek H_n$ is $\leq $ the number of graphs
in $\ek S_n$. Finally, it follows as in Lemma \ref{basicVnk}
that $H_n\ek H_n$ is a quotient of $H_n\ek \otimes_{H_{2k+1,n}} \ek H_n$,
where the latter has dimension $\leq \dim \Z[x] S_n\ek S_n$.
Hence the dimension of  $qB_n(r,q)\cong \bigoplus H_n\ek H_n$
is at most the one of the Brauer algebra. This proves the other
inequality.

To prove part (b), observe that we obtain a representation of
$\qB_n(r,q)$ with respect to the basis $(g_d)$ with coefficients in
$\Z[q^{\pm 1}, r^{\pm 1}, (r-1)/(q-1)]$. Specializing $r=q^N$,
these coefficients become elements of $\Z[q,q^{-1}]$ and we obtain
a representation $\pi$ of $\BnN$. As $\pi(g_d)1=g_d$, it follows
that the image has dimension at least $n!!$. The other inequality
follows as before from the fact that $(g_d)$ is a spanning set for
$\BnN$.

The proof of statement (c) follows from standard arguments.
Fix a basis $(g_d)$ and consider the left regular representation $\pi_l$
with respect to this basis. Then the discriminant
$\det(Tr(\pi_l(b_db_{d'})))$ is a polynomial in $r$ and $q$.
It specializes for $r=q^N$ and $q\to 1$ to the discriminant of $D_n(N)$,
which is known to be nonzero for $N>n$. This shows semisimplicity.
Similarly, the decomposition of a $\qB_n(r,q)$-module into simple
ones is already determined by the decomposition of any
specialization for $r$ and $q$, provided this specialized
algebra has the same decomposition into simple matrix algebras.

\begin{remark}\label{decompremark} If $V_{k,\nu}$ is a simple
$\qB_{k,\nu}$-module, we have the decomposition
$$V_{n,\la}\ \cong \ \bigoplus_{\mu} V_{n-1,\mu},$$
where $\mu$ runs through diagrams obtained by removing or also,
if $|\la|<n$, by adding a box to $\la$. This follows from
the restriction rule for the classical Brauer algebra,
essentially going back to Brauer (see also e.g. \cite{w2}).
If $|\la|=n$, this becomes
the restiction rule of modules of $S_n$ and $H_n$.
\end{remark}

\section{Markov trace}\label{SecMarkov}

\subsection{Definitions}
It will be convenient to slightly extend the ground rings.  So throughout this
section we will consider
the algebra $\BnN$ defined over the ring $\Z[q,q^{-1},[N]^{-1}]$,
and the algebra $\qB_n(r,q)$ defined over $\Z[q^{\pm 1}, r^{\pm 1},
((r-1)/(q-1))^{\pm 1}]$. For simplicity, we will only formulate the results
for $\BnN$; all the proofs will go though as well for $\qB_n(r,q)$.
 We  can now define the elements
$\eb=\frac{1}{[N]}e$ and $\ekb=\frac{1}{[N]^k}\ek$; for $B_n(r,q)$, we replace
$[N]^{-1}$ by $(1-q)/(1-r)$. Observe that
$\eb$ and $\ekb$ are idempotents with $\bar e_{(m)}\ekb=\ekb$ for $m\leq k$.
Recall that
a functional $\phi$ on an algebra $A$ has the trace property if
$\phi(ab)=\phi(ba)$ for all $a,b\in A$. It is well-known
that one can inductively define a trace functional $tr$ on $H_n$ by
$tr(1)=1$, and
$tr(g_{n-1}h)=\frac{q^N}{[N]}tr(h)$ for any $h\in H_{n-1}$.
Such a functional on
the Hecke algebras $H_n$ is called a {\it Markov trace}. It is compatible
with the obvious standard inclusion $H_{n-1}\subset H_n$.

\begin{lemma} \label{Markovprep1}
(a) There exists an isomorphism $\Psik$ between
$\ekb\BnN\ekb$ and $\qB_{n-2k}(N)$ such that $\Psik(\ekb g_i)=g_{i-2k}$
for $i>2k$ and $\Psik(\bar e_{(k+1)})=\eb$.

(b)  There exists a functional $\Phi_k:\BnN\to \qB_{n-2k}(N)$
uniquely defined by  $\Phi_k(h)=\Psik(\ekb h\ekb) $.
\end{lemma}

This lemma can be fairly easily checked using  Lemma
\ref{reppreparation} and the explicit basis
for $\qB_n$ in Theorem \ref{structuretheorem}.

\begin{lemma}\label{Markovdef} There exists a unique extension,
also denoted by $tr$ of the
Markov trace on $H_n$ to $\BnN$ which is defined via induction on $n$ by
$tr(a\ekb b)=tr(\ekb ba \ekb)=\frac{1}{[N]^{2k}}tr(\Phi_k(ba))$. This extension
also has the trace property  $tr(cd)=tr(dc)$ for all $c,d\in\BnN$.
\end{lemma}

$Proof.$  We will prove well-definedness and the trace property of
the functional $tr$ by induction on $n$. This is easy to check for
$n=1,2$, as the
algebras $\qB_1(N)$ and $\qB_2(N)$ are abelian. As to well-definedness
in general, we have to show that $tr(ac\ekb b)=tr(a\ekb cb)$ for all
$a,b\in H_n$ and $c\in H_{2k+1,n}$. This is equivalent to showing
$tr(\ekb ba \ekb c)=tr(c \ekb ba\ekb)$ by definition of $tr$.
But this follows from the trace property of $tr$ for $\qB_{n-2k}(N)$,
using the homomorphism $\Psi_k$.

Let us now prove the trace property for elements
$(a_1\bar e_{(k_1)}b_1)$ and $(a_2\bar e_{(k_2)}b_2)$, with
 $a_1,a_2,b_1,b_2\in H_n$.
Recall that we already know that
$tr(ab)=tr(ba)$ if $a,b\in H_n$.
Assuming $k_1\leq k_2$,  we can write
$$\bar e_{(k_1)}b_1a_2\bar e_{(k_2)}=
\sum_{j\geq k_2}a^{(j)}\bar e_{(j)}b^{(j)}$$
for suitable $a^{(j)}, b^{(j)}\in H_{2k_2+1,n}$. So we have
$$tr((a_1\bar e_{(k_1)}b_1)(a_2\bar e_{(k_2)}b_2)) =
\sum_{j\geq k_2}tr(a_1a^{(j)}\bar e_{(j)}b^{(j)}b_2)
=\sum_{j\geq k_2}tr(\bar e_{(j)}b^{(j)}b_2a_1a^{(j)}\bar e_{(j)})=$$
using $\bar e_{(k_2)}\bar e_{(j)}=\bar e_{(j)}$ for $j\geq k_2$
and $\bar e_{(k_2)}b^{(j)}=b^{(j)}\bar e_{(k_2)}$
$$=
\sum_{j\geq k_2}tr(\bar e_{(j)}b^{(j)}e_{(k_2)}b_2a_1e_{(k_2)}a^{(j)}\bar e_{(j)})
=\sum_{j\geq k_2}tr(e_{(k_2)}b_2a_1e_{(k_2)}a^{(j)}\bar e_{(j)}b^{(j)})=$$
$$
= tr((\bar e_{(k_2)}b_2a_1\bar e_{(k_2)})(\bar e_{(k_2)}b_1a_2\bar e_{(k_2)}))
= tr((\bar e_{(k_2)}b_1a_2\bar e_{(k_2)})(\bar e_{(k_2)}b_2a_1\bar e_{(k_2)})),$$
where we used the induction assumption
for elements in $\bar e_{(k_2)}\BnN\bar e_{(k_2)}\cong \qB_{n-2k_2}(N)$.
Equality with
$tr((a_2\bar e_{(k_2)}b_2)(a_1\bar e_{(k_1)}b_1))$
is now shown  by the same calculations as above.
Checking the trace property for elements $a\in H_n$ and $a_2\ek b_2$
goes similarly and is easier. The lemma is proved.

\subsection{Markov Property: Preparations} The goal is to prove an
analog of the Markov property for the extension of $tr$ to $\BnN$.
We will need the following technical lemmas:

\begin{lemma}\label{Markovprep2}
(a) If $j_1<i_1<n-1$,
$e\gm_{2,n-1}\gm_{1,j_1}g_n^{-1}\gp_{i_1,1}\gp_{n-1,2}e=
\gp_{i_1+1,3}g_{n,4}e_{(2)}\gm_{4,n}\gm_{3,j_1+2}$.

(b) If $i_1<j_1<n-1$, then
$e\gm_{2,n-1}\gm_{1,j_1}g_n^{-1}\gp_{i_1,1}\gp_{n-1,2}e=
\gp_{i_1+2,3}g_{n,4}e_{(2)}\gm_{4,n}\gm_{3,j_1+1}$.

(c) If $a,b\in H_n$, then $tr(a\gm_{n,2}\eb g_{2,n}b)=tr(ab)tr(\eb)$.
\end{lemma}

$Proof.$ Using Lemma \ref{braidrelations} and Eq \ref{Heckerel1},
we see that the left hand side
of statement (a) is equal to
\begin{align}\notag
e\gp_{i_1+1,3}\gm_{2,n}\gp_{n-1,2}e\gm_{3,j_1+2}\  & =\
\gp_{i_1+1,3}eg_2^{-1}g_1\gm_{3,n}\gp_{n-1,2}e\gm_{3,j_1+2} =\\
& =\ \gp_{i_1+1,3}\gp_{n,4}eg_2^{-1}g_1g_3^{-1}g_2e\gm_{4,n}\gm_{3,j_1+2},
\end{align}
which is equal to the right hand side of the statement. Statement (b) is proved
similarly.

For statement (c), one observes that any element $h\in H_{2,n}$ can be
written as a linear combination of elements in $H_{3,n}$ and
elements of the form $h_1g_2h_2$, with $h_1,h_2\in H_{3,n}$.
Then we have
$tr (\eb h_1g_2h_2\eb) = tr(h_1\eb g_2\eb h_2)= tr(g_2) tr(h_1\eb h_2)
= tr(\eb)tr(h_1g_2h_2),$
using relation $(E2)$ and the Markov property of $tr$ for Hecke algebras.
It follows that
$tr(\eb h \eb)=tr(h)tr(\eb)$ for any $h\in H_{2,n}$.
By  Lemma \ref{braidrelations}  the map $h\in H_n\mapsto
g_{1,n}h\gm_{n,1}\in H_{2,n+1}$ defines a trace-preserving homomorphism
from  $H_n$ onto $H_{2,n+1}$.
Claim (c) follows from this and the trace property.

\begin{lemma}\label{Markove}
Let $a,b\in H_n$. Then $tr(aebg_n^{-1})=tr(g_n^{-1})tr(aeb)$.
\end{lemma}

$Proof.$ We are going to prove the theorem by induction on $n$, with $n=1$ and $n=2$
easy to check. We will also need the fact that
$eH_ne\subset eH_{3,n}+H_{3,n}\e_{(2)}H_{3,n}$. Indeed
this can be checked easily using the fact that $H_n$ is the span
of elements of the form $g^-_{j_1,1}g_{j_2,2}h$ with $h\in H_{3,n}$.
Hence if the claim holds for $n-2$, then we also have $tr(g_n^{-1}ehe)
=tr(g_n^{-1})tr(ehe)$ by using the definition of $tr$ and induction assumption.

To prove the claim, let us write
$a=g_{i_1,1}g_{i_2,2}a''$ and $b=b''\gm_{2,j_2}\gm_{1,j_1}$,
where $a'',b''\in H_{3,n}$. We first observe that the claim
follows if both $i_1,i_2<n-1$.
Indeed, we have
\begin{align}
tr(a\eb bg_n^{-1}) &\ =\ & tr(g_{i_1,1}g_{i_2,2}g_n^{-1}\eb a''b)  &\ =\ &
tr(g_n^{-1}\eb a''bg_{i_1,1}g_{i_2,2}e)=  \notag \\ \notag
 &\ =\ &
tr(g_n^{-1})tr(\eb a''bg_{i_1,1}g_{i_2,2}\eb) &\ =\ &
   tr(g_n^{-1})tr(a\eb),  \\ \notag
\end{align}
where we used the argument of the first paragraph for the beginning of
the second line.
Similarly, one shows the claim if both $j_1,j_2<n-1$. Hence we can
assume that at least one of $i_1$ or $i_2$ is equal to $n-1$.
But as
$g_{n-1,1}g_{i_2,2}e=g_{i_2-1,1}g_{n-1,1}e=qg_{i_2-1,1}g_{n-1,2}e$,
we can assume that $i_2=n-1$ and $i_1<n-1$. One similarly shows
that we can assume $j_2=n-1$ and $j_1<n-1$.
Using Lemma \ref{Markovprep2} and the
isomorphism $\eb\qB_n\eb\cong \qB_{n-2}$, we can calculate for the case
$j_1<i_1$ that

\begin{align}
tr(a\eb bg_n^{-1})&\ =
\ & tr(b''\eb\gm_{2,n-1}\gm_{1,j_1}g_n^{-1}g_{i_1,1}g_{n-1,2}\eb a'')\ =\  &
 tr(b''g_{i_1+1,3}g_{n,4}\bar e_{(2)}\gm_{4,n}\gm_{3,j_1+2}a'')\ =\ \notag  \\ \notag
&\ = \ &
tr(\bar e_{(2)}\gm_{4,n}(\gm_{3,j_1+2}a''b''g_{i_1+1,3})g_{n,4}\bar e_{(2)})
\ = \ & tr(\gm_{3,j_1+2}a''b''g_{i_1+1,3})tr(\bar e_{(2)}).\notag
\end{align}
It remains to calculate
$tr(a\eb b)$. We get

\begin{align}
tr(a\eb b)&\ =
\ & tr(b''\eb\gm_{2,n-1}\gm_{1,j_1}g_{i_1,1}g_{n-1,2}\eb a'')& \ =\
 tr(b''g_{i_1+1,3}\eb \gm_{2,n-1}g_1 g_{n-1,2}\eb \gm_{3,j_1+2}a'')\ =\ \notag  \\
&\ = \ &
tr(b''g_{i_1+1,3}\eb g_2^{-1}g_1g_2\eb \gm_{3,j_1+2}a'' )
& \ = \  tr(g_2) tr(b''\gm_{3,j_1+2}\eb \gm_{3,j_1+2}a'')\ =\  \notag \\
&\ = \ &  tr(g_2)tr(b''\gm_{3,j_1+2}\eb \gm_{3,j_1+2}a'').&  \notag \\ \notag
\end{align}
The claim now follows from this and the fact that $tr(\bar e_{(2)})=
tr(g_2)(tr(g_n^{-1})tr(\eb))$. The case $i_1>j_1$ goes similarly, and $i_1=j_1$
is easy.

\subsection{Proof of Markov property}

\begin{theorem}\label{Markovprop}
The functional $tr$ satisfies $tr(cg_n)=tr(c)tr(g_n)$ for all
$c\in\BnN$.
\end{theorem}

$Proof.$ Observe that that the claim follows for $c\in H_n$ by
definition of $tr$, and for $c\in H_neH_n$ by Lemma \ref{Markove}.
We will prove the general claim by induction on $n$. It is trivially
true for $n=1$. If $n=2$, we have $tr(g_1g_2)=tr(g_1)tr(g_2)$ by
definition of $tr$, and $tr(\bar e g_2)=tr(\bar e g_2\bar e)=
\frac{q^n}{[N]}tr(\bar e)= tr(g_2)tr(\bar e)$ by relation $(E2)$.

Assuming that the claim holds for $n-1$ and $n-2$, we also have
$tr(\eb c\eb g_n)=tr(\eb c \eb)tr(g_n)$ for any $c\in H_n$,
using the isomorphism between $\eb \qB_{n+1}\eb$ and $\qB_{n-1}$,
see Lemma \ref{Markovprep1}.  The induction step in our proof will
depend on this observation.

Recall that any $b\in H_n$
can be written as $b=g_{i_{n-1},n-1}b'$ with
$b'\in H_{n-1}$ and $1\leq i_n\leq n$; here $g_{n,n-1}$ stands for 1,
i.e. $b=b'\in H_{n-1}$. But then we have
$$tr(a\ek b g_n)=tr(a\ek g_{i_n,n}b')=tr( b'a\ek g_{i_{n-1},n}).$$
One deduces that it suffices to show that $tr(a\ek g_{i_{n-1},n})=tr(g_n)
tr(a\ek g_{i_{n-1},n-1})$.
Now if $i_{n-1}>2$, $g_{i_{n-1},n-1}$ commutes with $\eb$ and we have
$$tr(a\ek g_{i_{n-1},n})=tr(a\ek  g_{i_{n-1},n-1}\eb g_n) =
tr(\eb a\ek  g_{i_{n-1},n-1}\eb)tr(g_n).$$
The claim now follows after verifying that the first factor in the last
expression is indeed equal to $tr(a\ek g_{i_{n-1},n-1})$.
As $\ek g_1=q\ek$, it only remains to consider the case $i_{n-1}=2$.
But then we have for $k\geq 2$, using Lemma \ref{Brrelations},(b)
that
$$tr(a\ekb g_{2,n})=tr((a\ekb g_2g_1)g_{4,n})=
tr((\eb g_2g_1a\ekb g_{4,n-1}\eb) g_n).$$
The claim now follows again by the argument mentioned at the beginning
of this proof.

\ignore{
For the induction step, observe that we have $tr(cg_n)=tr(c)tr(g_n)$
for any $c\in \bar e \BnN \bar e$; indeed, this follows from
Lemma \ref{Markovprep2} and the fact that $tr( \ekb h)=tr(\ekb)tr(h)
=\frac{1}{[N]^k}tr(h)$ for all $h\in H_{2k+1,n}\cong H_{n-2k}$.
We divide the proof into the following steps:
{\it Step 1:} We claim that it suffices to prove
$tr(a \ekb g_{1,n})=tr(a \ekb g_{1,n-1})tr(g_n)$. To see this,
first observe that in view of $\BnN=\bigoplus_k H_n\ek H_n$,
it suffices to prove the claim for $c=a\ek b$. Moreover, we can write
$b=g_{i_n,n-1}b'$ with $b'\in H_{n-1}$ and $1\leq i_n\leq n$, where
$i_n=n$ stands for the case  $b=b'\in H_{n-1}$. Observe that
$b'$ commutes with $g_n$ and
if $i_n>2$, then $eg_{i_n,n-1}=g_{i_n,n-1}e$. So we have
$$tr(a\ek b g_n)=tr(a\ek g_{i_n,n}b')=tr( g_{i_n,n}b'a\ek \bar e)
=tr( g_{i_n,n}\bar e b'a\ek);$$
We can now apply the induction assumption, as explained above
to see that the last expression is equal to
$$tr((\bar e b'a\ek g_{i_n,n-1})g_n)=tr(\bar e b'a\ek g_{i_n,n-1})tr(g_n)
=tr(a\ek b)tr(g_n),$$
for $i_n>2$. If $i_n=2$, we can use the fact that $\ek g_{2,n-1}=q^{-1}
\ek g_{1,n-1}$ to reduce the problem to the case $i_n=1$. Finally,
we can assume $b'=1$, as $tr(a\ek g_{1,n}b')=tr((b'a)\ek g_{1,n})$,
i.e. we replace $a$ by $b'a$.
{\it Step 2:} We now can show, for $k>1$ that
$$tr(a\ekb g_{2,n})=tr(a\ekb g_2g_1g_{4,n})=
tr(\eb g_2g_1a\ekb g_{4,n-1}\eb g_n).$$
The claim  follows again from the induction assumption.
}

\subsection{Weights} It is well-known that any trace functional on a
full $m\times m$ matrix algebra is equal to the usual trace, i.e. the
sum of the diagonal elements, up to a scalar multiple. Hence any trace
functional on a direct sum of full matrix algebras is completely
determined as soon as one knows this multiple for each summand;
these multiples are called the weights of the trace.
The weights for the Markov trace on the Hecke algebra $H_n$
for $tr(g_i)=r(q-1)/(r-1)$ and $\la$ a Young diagrams with $n$ boxes
are given by (see \cite{WHe})
\begin{align}\label{Heckeweights1}
\tilde\omega_\la\ & =\ q^{c_1(\la)}\ (\frac{q-1}{r-1})^n\
\prod_{(i,j)\in\la} \frac{rq^{i-j}-1}{q^{h(i,j)}-1}\ =\ \\ \notag
&=\  \frac{q^{c_2(\la)}}{[N]^n}\prod_{1\leq i<j\leq N}
\frac{[\la_i-\la_j+j-i]}{[i-j]}.
\end{align}
Here $c_1(\la)$ and $c_2(\la)$ are determined such that the formulas
remain invariant under the simultaneous
substitutions $r\mapsto r^{-1}$ and $q \mapsto q^{-1}$,
and equality with the second expression holds for $r=q^N$,
for Young diagrams with at most $N$ rows.
Moreover, $(i,j)$ denotes row and column of a box in the Young diagram $\la$,
$h(i,j)$ is the length of the hook in $\la$ with corner at $(i,j)$
given by
\begin{equation}
h(i,j)=\lambda_i-i+\lambda_j' -j+1,
\end{equation}
where $\la_i$ and $\la_j'$ denote the number of boxes in the $i$-th
row and $j$-th column of $\la$.  For more details, see e.g. \cite{Mac}.
Moreover, if $r=q^{N}$, we also have
\begin{equation}\label{GLchar}
[N]^n\tilde\om_\la=\chi_\la^{Gl(N)}(1,q,\ ...,\ q^{N-1}),
\end{equation}
where the right hand side is the character of an element of $Gl(N)$
with the indicated eigenvalues in the simple representation labeled by
$\la$.
We shall now similarly appeal to the character formulas
of orthogonal groups to calculate the weights of $tr$
for the algebras $\BnN$. We will need the following quantities
for a given Young diagram $\la$
\begin{equation}
d(i,j)\ = \begin{cases}
\lambda_i+\lambda_j-i-j & \text{if $i\leq j$,}\\
-\lambda_i'-\lambda_j'+i+j-2 &\text{if $i> j$.}
\end{cases}
\end{equation}

\begin{theorem}\label{Brweights}
The weights of the Markov trace $tr$ for $Br_n(r,q)$ are given by
\begin{align}\label{Heckeweights1}
\omega_{\la ,n}\ & =\ q^{c_3(\la)}\ (\frac{q-1}{r-1})^n\
\prod_{(i,j)\in\la} \frac{rq^{d(i,j)}-1}{q^{h(i,j)}-1}, \notag
\end{align}
where $\la$ runs through all the Young diagrams with $n, n-2, n-4, ... $
boxes, and $c_3(\la)$ is determined such that the formula is invariant
under the substitution $q\mapsto q^{-1}$.
\end{theorem}

$Proof.$
Recall that  the generic structures of $H_n$ and $Br_n(r,q)$ coincide with
the ones of the group algebra of the symmetric group and of the Brauer
algebra. Moreover, these isomorphisms are compatible  with the
inclusions. We have faithful representations of $S_n$ and $D_n(N)$
on $V^{\otimes n}$ if $N=\dim V >n$, where a minimal idempotent of
$\C S_n$ projects onto an irreducible representation of $Gl(N)$ in
$V^{\otimes n}$
and a minimal idempotent of $D_n(N)$ projects onto an irreducible
representation of $O(N)$.
Hence it follows
\begin{equation}\label{triangular}
\tilde\omega_\la\ =\ \sum_\mu b^\la_\mu \omega_{\mu,n},
\end{equation}
where $b^\la_\mu$ is the multiplicity of the irreducible $O(N)$-module
labeled by $\mu$ in the irreducible $Gl(N)$-module labeled by $\la$.
Moreover, we have $b^\la_\la=1$ and $b^\la_\mu\neq 0$ for $\mu\neq \la$
only if $\mu$ has fewer boxes than $\la$. Hence Eq.
\ref{triangular} gives us a triangular system of equations from which
we can calculate $\omega_\la$ for all $\la$s.
As
$$[N]^n\tilde\omega_\la = \chi_\la^{Gl(N)}(1,q,\ ...,\ q^{N-1})
= q^{n(N-1)/2}\chi_\la^{Gl(N)}(q^{(1-N)/2},q^{(3-N)/2},\ ...,\ q^{(N-1)/2})$$
for $r=q^N$,
we obtain the solution
$$\omega_\la =  \frac{1}{[N]^n}q^{n(N-1)/2}
\chi_\la^{O(N)}(q^{(N-1)/2},q^{(N-3)/2},\ ...,\ q^{(1-N)/2})
\quad {\rm if\ } r=q^N.$$
If $N$ is odd and sufficiently large,
the character on the right hand side is what is called the principal
character for type $B_{(N-1)/2}$ in \cite{Ko}.
It is shown in that paper that
$$\chi_\la^{O(N)}(q^{(1-N)/2},q^{(3-N)/2},\ ...,\ q^{(N-1)/2})
= q^{c_4(\la)}\ \prod_{(i,j)\in\la} \frac{q^{N+d(i,j)}-1}{q^{h(i,j)}-1},$$
with $c_4(\la)$ again chosen such that the formula is invariant under
the substitution $q\mapsto q^{-1}$.
Substituting $r=q^N$ in the numerators, we obtain the desired expression
for the weights. As these equalities hold for $r$ equal to any sufficiently
large odd power of $q$, they must hold true in general for rational
functions in $q$ and $r$.

\begin{remark} Contrary to the statement in \cite{Ko}, the principal
characters for type $B_n$ (and also for other types) do not coincide
with the $q$-dimensions of the corresponding quantum group
(the computations in the paper are correct, though).
The corresponding two-variable polynomials for these $q$-dimensions
have been calculated in \cite{w1} as $Q_\la(r,q)$ in connection with
another $q$-deformation of Brauer's centralizer algebra and lead
to different weights than
 the $\omega_{\la,n}$ in this paper.
\end{remark}

\subsection{Special values} The formulas for the weights of the Markov
trace are valid for the generic case, i.e. when $r$ and $q$ are viewed
as variables over a ring of rational functions. In this case, our algebras
are semisimple. These formulas will also hold if we define the algebras
$\qB_n$ over, say, the complex numbers, for any values of $r$ and $q$
for which $\qB_n(r,q)$ will have the same decomposition into a direct
sum of simple matrix rings as in the generic case. We shall use the weights
of the trace to determine these values, and also to determine special
semisimple quotients for certain cases when the algebras are not
semisimple.

We define special finite sets $\LNell$ of Young diagrams
for integers $N$ and $\ell$ satisfying $1<|N|<\ell$. These will be
related to algebras $\qB_n(r,q)$ where $r=q^N$ and $q=\xi$ is a primitive
$\ell$-th root of unity.

\begin{definition} Fix  integers $N$ and $\ell$ satisfying $1<|N|<\ell$.
The set $\LNell$ consists of all Young diagrams
$\la$ with $\la_i$ boxes in the $i$-th row and $\la_j'$ boxes in
the $j$-th column which satisfy

(a) $\la_1'+\la_2'\leq N$ and $\la_1\leq (\ell-N)/2$ if
$N>0$ and $\ell-N$ even,

(b)  $\la_1'+\la_2'\leq N$ and $\la_1+\la_2\leq \ell-N$ if
$N>0$ and $\ell-N$ odd,

(c) $\la_1\leq |N|/2$ and $\la_1'+\la_2'\leq \ell-|N|$ if
$N<0$ is even,

(d) $\la_1+\la_2\leq |N|$ and $\la'_1+\la'_2\leq \ell-|N|$ if
$N<0$ is odd.

In each of these cases, we call a Young diagram a boundary diagram
of $\LNell$ if it satisfies one inequality of the definition,
but misses the other one
by 1 (e.g. in case (a) if $\la_1'+\la_2'\leq N$ and $\la_1=1+(\ell-N)/2$.
We denote by $\LNellb$ the union of $\LNell$ with its boundary diagrams.
\end{definition}

\begin{proposition}\label{specialval} (a) The weights $\om_{\la ,n}=
\om_{\la ,n}(\xi^N,\xi)$
are nonzero and well-defined for any primitive $\ell$-th root of unity $\xi$.

(b) If $\xi$ is a primitive $\ell$-th root of unity, then
$\om_{\la ,n}(\xi^N,\xi)\neq 0$ for $\la\in\LNell$, and
$\om_{\la ,n}(\xi^N,\xi)=0$ for any boundary diagram of $\LNell$.
\end{proposition}

$Proof.$ The statements can be easily checked using the explicit product
form of the formulas for $\om_{\la,n}$.

\begin{lemma}\label{Heckesemisimple} Let $\xi$ be a primitive $\ell$-th root
of unity and let $|N|\geq 2$.

(a) Every Specht module $S_\la$ of the Hecke algebra $H_m(\xi)$ labeled by a
Young diagram $\la$  in $\LNellb$ with $m$ boxes is simple.

(b) If $V$ is an $H_n(\xi)$-module which decomposes as an $H_{n-1}(\xi)$-module
into a direct sum of simple modules labeled by Young diagrams in $\LNell$ with
$n-1$ boxes, then $V$ is also semisimple as an $H_n(\xi)$-module, with its
simple components labeled by Young diagrams in $\LNellb$.
\end{lemma}

$Proof.$  It follows from the Nakayama Conjecture for Hecke algebras
(a theorem proved in \cite{DJ}) that any Specht module is simple if
it is labeled by a
Young diagram $\la$ for which $\la_1+\la_1'<\ell+1$.
Moreover, the corresponding central idempotent $z_\la$ is well-defined
for a primitive $\ell$-th root of unity. This can also be
easily checked using the explicit representations e.g. in \cite{WHe}.
Statement (a) can now be fairly easily checked using this criterion.

To prove statement (b), let $z_{(n-1)}=\sum_\mu z_\mu$, with $\mu$ in
$\Lambda(N-1,\ell)$, and let $z_{(n)}=\sum_\la z_\la$, with $\la\in\LNellb$.
It follows from the well-known restriction rule for simple Hecke algebra
modules in the semisimple case that $z_{(n)}z_{(n-1)}=z_{(n-1)}$.
Hence $V=z_{(n-1)}V=(z_{(n)}z_{(n-1)})V=z_{(n)}V$, also for $\xi$ a
primitive $\ell$-th root of unity. This proves part (b).

\section{Semisimplicity}\label{SecSemisimp}

We now view our algebras $\qB_n(\r,\xi)$ defined over a field of characteristic 0.
We determine for which values of the parameters $r=\r$ and $q=\xi$
in the chosen field our algebras
will be semisimple. This follows the same patterns as in \cite{w1} and
\cite{w2}, using Jones' basic construction and our formulas for
the weights of the trace from the previous section. The only new
complications come from the fact that we will
not be able to use the standard embeddings $\qB_n\subset\qB_{n+1}$.
We will often just write $\qB_n$ instead of $\qB_n(r,\xi)$, assuming
$\r$ and $\xi$ to be fixed.

\subsection{Jones' construction}\label{basicconstruction}
Let $A\subset B\subset C$ be finite
dimensional algebras. Moreover, let $tr$ be a trace functional on $B$
such that the induced bilinear form $(b_1,b_2)=tr(b_1b_2)$ is
nondegenerate for $B$, and also for its restriction on $A$.
We can then define a conditional expectation $E_A:B\to A$ uniquely
determined by
\begin{equation}\label{condexp}
(E_A(b),a)\ =\ (b,a)\quad {\rm for\ all\ }a\in A.
\end{equation}
Moreover, we assume that there exists an idempotent $p$ in $C$
satisfying the following conditions
\vskip .2cm
(a)\ $pa=ap$ for all $a\in A$, and the map $a\in A\mapsto ap$ is a
monomorphism,

(b)\ $pbp=E_A(B)p$ for all $b\in B$.
\vskip .2cm
Under these conditions we have the following results, going back
to Jones' basic construction (see \cite{w2}, Lemma 1.1 or \cite{w1},
Theorem 1.1):

\begin{proposition}\label{basicconstruct}
The ideal $\langle p \rangle$ in the algebra
generated by $B$ and $p$ is isomorphic to the commutant
$\End_A(B)$ of $A$, acting via right multiplication on $B$. In particular,
if $A$ is semisimple, so is  $\langle p \rangle$. Moreover, the ideal
$\langle p \rangle$ is spanned by elements of the form $b_1pb_2$, with
$b_1,b_2\in B$.
\end{proposition}

\subsection{Embeddings} We define the embeddings
$i_1, i'_1:\qB_{n-1}\to \qB_n$ by
$i_1(b)=g_{1,n-1}bg_{1,n-1}^{-1}$ and
$i'_1(b)=\gm_{1,n-1}b(\gm_{1,n-1})^{-1}$
for $b\in \qB_{n-1}$. Moreover,
we also define $i_2,i'_2:\qB_{n-2}\to \qB_n$ by
$i_2(b)=i_1(i'_1(b))=\gm_{2,n-1}i_1(b)(\gm_{2,n-1})^{-1}$,
for $b\in \qB_{n-2}$, and $i'_2(b)=i'_1(i_1(b))$.
Observe that we have $i_1(g_j)=g_{j+1}=i'_1(g_j)$ for $1\leq j<n-1$.
Then we have the following easy lemma:

\begin{lemma}\label{iprep} (a) With notations above we have
that  $i_2(\qB_{n-2})$ commutes with $e$ and the map
$b\in i_2(\qB_{n-2}) \mapsto b\bar e$ defines an injective homomorphism.
The statement also holds for  $i'_2(\qB_{n-2})$ instead of
$i_2(\qB_{n-2})$.

(b) Assume that $\qB_{n-1}$ is spanned by
elements of the form $b_1\chi b_2$ and that
$tr(b_1\chi b_2)=tr(\chi)tr(b_1b_2)$, where $\chi\in\{ 1, g_1, e\}$,
and $b_1,b_2\in i'_1(\qB_{n-2})$. Then we have
$\eb(i_1(b_1\chi b_2))\eb=tr(\chi)i_1(b_1b_2)\eb$ and also
$tr(c_1\chi c_2)=tr(\chi)tr(c_1c_2)$ for $c_1,c_2\in i_1(\qB_{n-1})$.

(c) Under the assumptions and notations  of (b),
we have $E_{i_1'(\qB_{n-2})}(b_1\chi b_2)=
tr(\chi)b_1b_2$, assuming that $tr$ induces nondegenerate bilinear
forms on $\qB_{n-1}$ and $\qB_{n-2}$.
\end{lemma}

$Proof.$ It follows from Lemma \ref{braidrelations} that
 $i_1(g_j)=i'(g_j)=g_{j+1}$ for $j<n-1$ and
$i_2(g_j)=i'_2(g_j)=g_{j+2}$ for $j<n-2$.
If we define $e_2=i_1(e)$, and $e_3=i_2(e)$, then it follows
from our relations that
$$ee_3=eg_2^{-1}g_3{-1}g_1g_2eg_2^{-1}g_1^{-1}g_3g_2=
e_{(2)}g_2^{-1}g_1^{-1}g_3g_2=e_{(2)}.$$
One similarly checks that $e_3e=e_{(2)}$. This, together with
the relation $eg_j=g_je$ for $j>2$ shows that $e$  commutes with
$i_2(\qB_{n-2})=A$. Hence the map $b\in \qB_{n-2}\mapsto
\eb i_2(b)$ is an algebra homomorphism. One checks easily
at the generators that it is the inverse of the isomorphism $\Psi_1$,
as defined in Lemma \ref{Markovprep1}.
The same proof goes through if we replace $i_2$ by $i'_2$.
This proves part (a).

For part (b), observe that $e_2=i_1(e)=g_1g_2eg_2^{-1}g_1^{-1}$.
If $\Delta_k=g_{1,k-1}g_{1,k-2}\ ...\ g_1$, then
$\Delta_n^{-1}i_1(e)\Delta_n=\Delta_{n-1}^{-1}e\Delta_{n-1}\in\qB_{n-1}$,
and $\Delta_n^{-1}g_i\Delta_n=g_{n-i}$.
One deduces from this that $\Delta_n^{-1}i_1(\qB_{n-1})\Delta_n=\qB_{n-1}$.
But then, if $b\in i_1(\qB_{n-1})$, we have
$$tr(g_1b)=tr(\Delta_n^{-1}g_1b\Delta_n) =
tr(g_{n-1}\Delta_n^{-1}b\Delta_n)=tr(g_1)tr(b),$$
using the trace property and Theorem \ref{Markovprop}.
Hence we only need  to prove the last statement of (b) for $\chi=e$,
or, equivalently, $\chi=\eb$.
By our assumptions, we can write $c_2c_1=i_1(b_1\psi b_2)$,
with $\psi\in \{ 1,e,g_1\}$ and $b_1,b_2\in i_1'(\qB_{n-1})$.
But then
$$tr(c_1\eb c_2)=tr(\eb i_1(b_1\psi b_2)\eb)=tr(\psi)tr(i_1(b_1)\eb i_1(b_2))=
tr(\eb)tr(b_1\psi b_2),$$
using our assumptions and already proven claims.
It only remains to prove claim (c), which follows from the definitions
and from  $tr(b_1\chi b_2c)= tr((tr(\chi)b_1b_2)c)$ for any $c\in
i_1'(\qB_{n-2})$.

\begin{theorem}\label{semisimplicity}
The algebra $\qB_n(\r,\xi)$ is semisimple if $\r\neq \xi^k$ for $|k|\leq n$
and if $\xi$ is not an $\ell$-th root of unity, $\ell \leq n$.
 In this case, it has the same decomposition into simple matrix
rings as the generic Brauer algebra, and the trace $tr$ is nondegenerate.
In particular, the assumptions in Lemma \ref{iprep} hold for all $n$.
\end{theorem}

$Proof.$
We will prove the claim by induction on $n$ together
with the spanning assumption in  Lemma \ref{iprep},(b),
with $n$ replaced by $n+1$ (i.e. when $b_1,b_2$ are in $i_1'(\qB_{n-1})$).
This, as well as the claim in the statement
is easy to check for $n=1$ and $n=2$.

By induction assumption, $tr$ is nondegenerate on $\qB_{n-1}$ and $\qB_n$.
Hence, by Lemma \ref{iprep}, all the assumptions for
 Prop. \ref{basicconstruct} are satisfied for $A=i_2(\qB_{n-1})$,
$B=i_1(\qB_{n})$ and $p=\eb$. Hence
the ideal $\langle e \rangle$ generated by $e$ in the algebra generated
by $i_1(\qB_n)$ and $e$ is isomorphic to $\End_{\qB_{n-1}}\qB_n$.
It is known from the generic Brauer algebra that the latter algebra
has dimension $(2n+1)!!-(n+1)!$; it is spanned by
all graphs which have at least one horizontal edge.
Using the basis $(g_d)$ of Theorem \ref{structuretheorem}, we
see that this ideal coincides with the ideal $I_{n+1}$ generated
by $e$ in $\qB_{n+1}$, and that it has zero intersection with $H_{n+1}$.
Now both $I_{n+1}$ and $H_{n+1}\cong
\qB_{n+1}/I_{n+1}$ are semisimple algebras with mutually nonisomorphic
simple modules (as $e$ acts nonzero on simple $I_{n+1}$-modules and zero
on simple $H_{n+1}$-modules). It follows that
 $\qB_{n+1}\cong I_{n+1}\oplus H_{n+1}$ as algebras.
Nondegeneracy of a trace on a semisimple algebra now can be checked
by just showing that its values on minimal idempotents are nonzero.
This follows from Theorem \ref{Brweights}.

Additionally, it follows from Prop. \ref{basicconstruct} and well-known
properties of the Hecke algebra $H_{n+1}$ that $\qB_{n+1}$ is spanned by
elements of the form $b_1\chi b_2$, with $b_1, b_2\in B=i_1(\qB_n)$
and $\chi\in \{ 1, e, g_1\}$. To prove the spanning assumption,
we observe that everything in this proof so
far would have worked as well for the inclusion $A'=i'_2(\qB_{n-1})\subset
B'=i'_1(\qB_{n})\subset \qB_{n+1}$. Hence $\qB_{n+1}$ is also spanned
by elements of the form  $b'_1\chi b'_2$, with $b'_1, b'_2\in B'=i'_1(\qB_n)$
and $\chi\in \{ 1, e, g_1\}$. This finishes the proof.

\begin{corollary}\label{semisimplecor} Let $Ann_n(\r,\xi)=\{ a\in \qB_n(\r,\xi), tr(ab)=0$ for all $b\in
\qB_n(\r,\xi)\}$ and let $\qBb_n(\r,\xi)=\qB_n(\r,\xi)/Ann_n(\r,\xi)$. Then
 $Ann_n(\r,\xi)\subset Ann_{n+1}(\r,\xi)$ for all $n$.
\end{corollary}

$Proof.$ Let $\Delta_{n+1}$ be as defined in the proof of Lemma
\ref{iprep}. We have seen in the proof of Theorem \ref{semisimplicity}
that $\qB_{n+1}$ is spanned by elements of the form $b_1\chi b_2$,
with $b_1,b_2\in i_1(\qB_n)$, and $\chi\in \{ 1,e,g_1\}$.
Conjugating this
by $\Delta_{n+1}$, we see that  $\qB_{n+1}(r,q)$ is also spanned by
elements of the form $c_1\psi c_2$, with $c_1,c_2\in \qB_n$ and
$\psi\in \{ 1,e_n=\Delta_{n+1}e\Delta_{n+1}^{-1},g_n\}$.
If $a\in Ann_n(\r,\xi)$, then we also have $tr(ac_1\chi c_2)=tr(\chi)tr(ac_1c_2)
=0$. Hence also $a\in Ann_{n+1}(\rho,\xi)$.

\begin{theorem} Let $\xi$ be a primitive $\ell$-th root of unity,
and let  $N$ be an integer satisfying $1<|N|<\ell$.
Then $\qBb_n(\xi^N,\xi)$ is semisimple for all $n\in\N$. Its simple
components are labeled by the Young diagrams in $\LNell$ with
$n,n-2, n-4, ...$ boxes, and the values of the Markov trace
for minimal idempotents in  $\qBb_n(\xi^N,\xi)$ are given by
the formulas in Theorem \ref{Brweights}. The restriction rule
from  $\qBb_n(\xi^N,\xi)$ to $\qBb_{n-1}(\xi^N,\xi)$ is as  in  Remark
\ref{decompremark}, where now
only diagrams in $\LNell$ are allowed.
\end{theorem}

$Proof.$ We will only write $\qB_n$ for $\qB_n(\xi^N,\xi)$ in this proof,
which  will be done by induction on $n$ similar to the one of
Theorem \ref{semisimplicity}. For $n=1$ and $n=2$,
the claim is easily checked. To prove the induction step $n\to n+1$,
we obtain from Corollary \ref{semisimplecor} that also $\qBb_{n+1}$
is semisimple, with the ideal
 $\langle e \rangle\cong \End_{\qBb_{n-1}}\qBb_{n}$, and
$\qBb_{n+1}\cong \langle e \rangle\oplus \bar H_{n+1}$,
where $\bar H_{n+1}$ is a quotient of the Hecke algebra $H_{n+1}$.
Moreover,
it is well-known in the setting of Section \ref{basicconstruction}
that we get minimal idempotents in $\End_AB$ in the form $pf$, where
$f$ is a minimal idempotent in $A$, acting from left on $A$.
Hence we get minimal idempotents in the ideal  $\langle e \rangle$
of the form $p_\la e$, where $p_\la$ is a minimal
idempotent in $i'_2(\qBb_{n-1,\la})\cong \qBb_{n-1,\la}$ with
$\la\in \LNell$ such that $n-1-|\la|$ is nonnegative and even.
We have $tr(p_\la e)=tr(e)tr(p_\la)$, as claimed.

It remains to determine the remaining simple components of $\bar H_{n+1}$.
By induction assumption and the restriction rules, see
Remark \ref{decompremark} and Lemma \ref{Heckesemisimple},
such a simple module must be isomorphic to a Specht module
labeled by a Young diagram $\la$  in  $\bar\Lambda(N+1,\ell)$.
So now it suffices to show that the trace of a
minimal idempotent in the corresponding simple component
is again given by $\om_{\la,n}$. This follows as soon as we can
find an explicit expression for a minimal idempotent
in $\qB_{N+1,\la}$ in terms of
basis elements for which the coefficients are rational functions
in $r$ and $q$ which are well-defined for our special values.
This can be done by using the path idempotent approach, as it
was done in \cite{RW}, as follows: Let $\mu$ be a diagram in
$\LNell$ obtained by removing a box from $\la$. It follows from
the restriction rule that the minimal idempotent
$p_\mu\in \qB_{n,\mu}$ can be written as a sum of
mutually commuting minimal idempotents $p_\nu\in\qB_{n+1,\nu}$
labeled by diagrams $\nu$ obtained by adding or subtracting a box
to/from $\la$. Now if $\nu$ has one box less than $\la$,
$\qB_{n+1,\nu}$ is in the basic construction part of $\qB_{n+1}$,
and hence $p_\nu$ can be obtained via formulas in \cite{RW}; see
\cite{RW}, Theorem 1.4 and our explicit formulas for the weights of the trace,
Theorem \ref{Brweights}.
In particular, they are well-defined at our given root of unity $q$.
Let $p_\mu'$ be the idempotent
obtained after  subtracting these idempotents $p_\nu$ from $p_\mu$.
We then obtain $p_\la$ as an eigenprojection from $p_\mu'g_np_\mu'$
using the formulas in \cite{WHe}, Cor. 2.3. This finishes the proof.

\begin{remark} 1. Essentially by the same method, semisimple quotients
were constructed
in \cite{w1} for another $q$-deformation of Brauer's centralizer algebra.
As the weights for the Markov traces for these two generalizations
of Brauer's algebras differ, we also obtain different quotients.
However, as in \cite{w1}, we will be able to construct new subfactors
of the hyperfinite II$_1$ factor from our algebras by exhibiting
a $C^*$-structure for certain quotients. This analysis will be
similar to the one in \cite{w1}, but the subfactors will be substantially
different. E.g. it is expected that for $N=2$ we would get the
Goodman-de la Harpe-Jones subfactors labeled by Dynkin graphs
$D_{2n}$, see \cite{GHJ}. This will be done in a future paper.

2. The semisimple quotients constructed in this paper are not
maximum in general. It is expected that the algebras in this paper
are cellular in the sense of \cite{GL}. It would be interesting to
determine their decomposition series.

3. It is possible to define a $q$-deformation of $U\so_N$ as a subalgebra
of $U_q\sl_N$, see \cite{L1},\cite{L}, \cite{IK}. It is not a sub Hopf algebra of $U_q\sl_N$
but a coideal algebra.  Hence its representations can be made
into a module category of $Rep(U_q\sl_N)$.
Taking the commutant of its action on $V^{\otimes n}$,
where $V$ is the vector representation, we obtain a $q$-deformation of
the Brauer algebra. This algebra was already studied in \cite{Mo}
(see remarks below) and is closely related to our algebras here.
In particular, as these coideal algebras were constructed
for a wide class of subalgebras of a semisimple Lie algebra, it might be
possible to generalize constructions of this paper in this more
general context. This would require more detailed studies of their
representation theory in the nonsemiple case.

4. {(\it Module categories)} It follows from the description via generators and relations that
the map $b\otimes g_i\mapsto bg_{i+m}$ defines
embeddings of $\qBb_m(\xi^N,\xi)\otimes\bar H_n(\xi)\subset \qBb_{n+m}(\xi^N,\xi)$,
with the algebras as defined in this section.
This should lead to the construction of a  module category of
the fusion tensor category of type $A_{|N|-1}$ of level $\ell-|N|$
(see e.g. \cite{Os}), with the objects being
idempotents of the algebras $\qBb_n(\xi^N,\xi)$. Here the fusion tensor
category would be defined via idempotents in the Hecke algebra quotients,
see e.g. \cite{Bl}. It appears that for $N=2$, we would obtain the module
tensor categories as in \cite{Os} given by Dynkin graphs $D_n$.
At least in this case, this category should also be realized
via bimodules of  von Neumann factors and subfactors
as mentioned in Remark 1.  Finally, we also mention that we obtain for each
set $\LNell$ a representation of the fusion ring of type $A_{|N|-1}$ of level
$\ell -|N|$ via matrices with nonnegative integer entries whose rows and columns
are labeled by the entries of $\LNell$. They describes the tensor product rules
of the model action. So our paper gives a rigorous derivation of at least some of the
NIMREP representations in e.g. \cite{GG} (see also the references in that paper).
This was one of the motivations for this paper.

5. It would be interesting to see whether our algebras have any topological
meaning. There  exist other algebras, motivated by topological considerations,
which contain Hecke algebras as unital subalgebras,
see \cite{Ju}, \cite{RH}. It is not clear at this point what the relation is
between these and our algebras, if any.

6. While putting on the finishing touches on this paper, the author noticed
the work \cite{Mo} by A.Molev. It deals with algebras acting on
tensor spaces which also are $q$-deformations of quotients
of Brauer's centralizer algebras. The structure analysis in \cite{Mo}
was considerably less detailed than in this paper, though.
It was conjectured in an earlier version of this
paper that those algebras should be related to the ones in this paper.
The author would like to thank A. Molev for informing him that
indeed the algebras in this paper do satisfy the relations of the ones in \cite{Mo}.

\end{remark}
\vskip .5cm
\bibliographystyle{plain}

\end{document}